\documentstyle[12pt]{article}
\def\Bbb R{{\rm \bf R}}
\def\proclaim#1{\vskip2mm{\bf #1}\em}
\def\endproclaim{\em \vskip2mm}
\def\tag#1{\eqno(#1)}
\def\gathered{\begin{array}{c}}
\def\endgathered{\end{array}}
\def\text{\mbox}

\begin{document}

\title {The enclosure method for
inverse obstacle scattering using a single electromagnetic wave in time domain}
\author{Masaru IKEHATA\footnote{
Laboratory of Mathematics,
Institute of Engineering,
Hiroshima University,
Higashi-Hiroshima 739-8527, JAPAN}}
%\date{}
\maketitle

\begin{abstract}
In this paper, a time domain enclosure method for an inverse
obstacle scattering problem of electromagnetic wave is introduced.
The wave as a solution of Maxwell's equations is generated by an applied volumetric current having an {\it orientation} and
supported outside an unknown obstacle and observed
on the same support over a finite time interval.  It is assumed that
the obstacle is a perfect conductor. Two types of analytical
formulae which employ a {\it single} observed wave
and explicitly contain information about the geometry of the
obstacle are given.  In particular, an effect of the orientation of the current is catched in one of two formulae.
Two corollaries concerning with the detection of the points on the surface
of the obstacle nearest to the centre of the current support
and curvatures at the points are also given.

\noindent
AMS: 35R30, 35L50, 35Q61, 78A46, 78M35

\end{abstract}

%\tableofcontents

\section{Introduction}

In this paper, we consider an inverse obstacle scattering problem
of a wave whose governing equation is given by Maxwell's equations.
The wave is generated by a {\it source} at $t=0$ which is {\it not
far a way} from an unknown {\it obstacle}, and we observe a {\it
single} reflected wave from the obstacle over a {\it finite time
interval} at the same place as the source.
The inverse obstacle scattering problem is to: extract information
about the geometry of the obstacle from the observed wave.  This
is a proto-type of so-called inverse obstacle problem \cite{IS} and the
solution may have possible applications to radar imaging.  Since
we consider the data over a finite time interval and thus,
this is a {\it time domain} inverse problem. Our main interest is
to find an analytical method or formula that extracts the geometry of the
obstacle from the data by using the governing equation of the
wave.

Let us describe the mathematical formulation of the problem.
Let $D$ be a nonempty bounded open subset of $\Bbb R^3$ with $C^2$-boundary such that
$\Bbb R^3\setminus\overline D$ is connected. $\mbox{\boldmath $\nu$}$ denotes the unit normal to
$\partial D$, oriented towards the exterior of $D$.

Let $0<T<\infty$.
We denote by $\mbox{\boldmath $E$}$ and $\mbox{\boldmath $H$}$ the electric field and the magnetic field, respectively.
$\epsilon$ denotes the electric permittivity and $\mu$ the magnetic permeability assumed to be positive
constants.

We assume that $\mbox{\boldmath $E$}$ and $\mbox{\boldmath $H$}$
are induced only by the current density $\mbox{\boldmath $J$}$ at $t=0$ and that the obstacle is a {\it perfect conductor}.
It is well known that the governing equations of $\mbox{\boldmath $E$}$ and $\mbox{\boldmath $H$}$
take the form
$$
\left\{
\begin{array}{ll}
\displaystyle
\epsilon\frac{\partial\mbox{\boldmath $E$}}{\partial t}
-\nabla\times\mbox{\boldmath $H$}=\mbox{\boldmath $J$}
& \text{in}\,(\Bbb R^3\setminus\overline D)\times\,]0,\,T[,\\
\\
\displaystyle
\mu\frac{\partial\mbox{\boldmath $H$}}{\partial t}
+\nabla\times\mbox{\boldmath $E$}=\mbox{\boldmath $0$}
& \text{in}\,(\Bbb R^3\setminus\overline D)\times\,]0,\,T[,\\
\\
\displaystyle
\mbox{\boldmath $\nu$}\times\mbox{\boldmath $E$}=\mbox{\boldmath $0$} &
\text{on}\,\partial D\times\,]0,\,T[,\\
\\
\displaystyle
\mbox{\boldmath $E$}\vert_{t=0}=\mbox{\boldmath $0$} & \text{in}\,\Bbb R^3\setminus\overline D,
\\
\\
\displaystyle
\mbox{\boldmath $H$}\vert_{t=0}=\mbox{\boldmath $0$}
& \text{in}\,\Bbb R^3\setminus\overline D.
\end{array}
\right.
\tag {1.1}
$$

Now let us describe our problem.
Fix a large (to be determined later) $T<\infty$.
Let $B$ be the open ball centred at a point $p$ with {\it very small} radius $\eta$
and satisfy $\overline B\cap\overline D=\emptyset$.
There are several choices of the current density $\mbox{\boldmath $J$}$
as a model of the antenna (\cite{B, CB}).
In this paper, we assume that $\mbox{\boldmath $J$}$
takes the form
$$\displaystyle
\mbox{\boldmath $J$}(x,t)
=f(t)\chi_B(x)\mbox{\boldmath $a$},
\tag {1.2}
$$
where $\mbox{\boldmath $a$}\not=\mbox{\boldmath $0$}$ is a constant unit vector,
$\chi_B$ denotes the characteristic function of $B$ and $f\in
H^1(0,\,T)$ with $f(0)=0$. Note that $\chi_B(x)$ has discontinuity
across the sphere $\partial B$.

{\bf Problem.}  Generate $\mbox{\boldmath $E$}$ and $\mbox{\boldmath $H$}$ by $\mbox{\boldmath $J$}$ and
observe $\mbox{\boldmath $E$}$ on $B$ over time interval $]0,\,T[$.
Extract information about the geometry of $D$ from the observed data.

This may be the simplified model of the case when the reflected
wave is observed at the same place where the source is located.
Note that we consider the pair $(\mbox{\boldmath $E$},
\mbox{\boldmath $H$})$ is the solution of (1.1) in the sense as
described on pages 433-435 in \cite{DL5} which is based on Stone's
theorem. In particular, we make use of the fact that
$(\mbox{\boldmath $E$}, \mbox{\boldmath $H$})$ belongs to
$C^1([0,\,T], L^2(\Bbb R^3\setminus\overline D)^3\times L^2(\Bbb
R^3\setminus\overline D)^3)$ with $(\nabla\times\mbox{\boldmath
$E$}(t),\nabla\times\mbox{\boldmath $H$}(t))\in L^2(\Bbb
R^3\setminus\overline D)^3\times L^2(\Bbb R^3\setminus\overline
D)^3)$ and $\mbox{\boldmath $E$}(t)\times\mbox{\boldmath
$\nu$}\vert_{\partial D}=\mbox{\boldmath $0$}$ for all
$t\in\,[0,\,T]$.

As far as the author knows there is no result for the problem mentioned above.
The point is: the data is taken over a {\it finite time interval} and
only a {\it single} (reflected) wave is employed.

In this paper, we employ the {\it enclosure method} for this
problem. The origin goes back to a method developed for an inverse
boundary value problem in two dimensions for the Laplace equation
\cite{I1}.

The method consists of two tools:

$\bullet$  a special solution  $v$ of an elliptic
partial differential equation which depends on a large parameter
$\tau>0$ and is independent of unknown obstacles.

$\bullet$  a so-called indicator function of independent variable $\tau$
constructed by using observation data and $v$ above.

Studying the asymptotic behaviour of the indicator function as
$\tau\longrightarrow\infty$ yields some information about the
location and shape of unknown discontinuity.

We have already some applications to inverse obstacle scattering problems
whose governing equation is given by the classical wave equation
in three-space dimensions \cite{IE, IE2, IE3, IE4, IE5}. The method
enables us to extract information about the geometry of unknown
obstacle from a {\it single} reflected wave over a finite time
interval. However, the governing equation therein is a single
partial differential equation and it is not clear that the method
can cover also the very important case when the governing equation
consists of a system of partial differential equations.

In the following subsection we describe our solution to Problem.

\subsection{Statement of the results}

We denote by $H(\text{curl},\,\Bbb R^3)$ the set of all vector valued-functions $U\in L^2(\Bbb R^3)^3$
such that $\nabla\times\mbox{\boldmath $U$}\in L^2(\Bbb R^3)^3$.
It is a Hilbert space with norm
$$\displaystyle
\Vert U\Vert_{H(\text{curl},\,\Bbb R^3)}
=\sqrt{\Vert\mbox{\boldmath $U$}\Vert_{L^2(\Bbb R^3)^3}^2+
\Vert\nabla\times\mbox{\boldmath $U$}\Vert_{L^2(\Bbb R^3)^3}^2}
$$
and $C_0^{\infty}(\Bbb R^3)^3$ is dense in $H(\text{curl}\,,\Bbb R^3)$.

By the Lax-Milgram theorem, we know that given $\mbox{\boldmath $f$}(\,\cdot\,,\tau)\in L^2(\Bbb R^3)^3$
there exists a unique $\mbox{\boldmath $V$}\in H(\text{curl}\,,\Bbb R^3)$ such that, for all $\mbox{\boldmath $\Psi$}\in H(\text{curl}\,,\Bbb R^3)$
$$\displaystyle
\int_{\Bbb R^3}
\left(\frac{1}{\mu\epsilon}\nabla\times\mbox{\boldmath $V$}\cdot\nabla\times\mbox{\boldmath $\Psi$}+\tau^2\mbox{\boldmath $V$}\cdot\mbox{\boldmath $\Psi$}
\right)dx+\int_{\Bbb R^3}\mbox{\boldmath $f$}(x,\tau)\cdot\mbox{\boldmath $\Psi$}dx=0.
$$
We call this $\mbox{\boldmath $V$}$ the weak solution of
$$\begin{array}{c}
\displaystyle
\frac{1}{\epsilon\mu}\nabla\times\nabla\times\mbox{\boldmath $V$}+\tau^2\mbox{\boldmath $V$}
+\mbox{\boldmath $f$}(x,\tau)
=\mbox{\boldmath $0$}\,\,\text{in}\,\Bbb R^3.
\end{array}
\tag {1.3}
$$
In this paper, unless otherwise stated, $\mbox{\boldmath $f$}(\,\cdot\,,\tau)$ has the form
$$\displaystyle
\mbox{\boldmath $f$}(x,\tau)
=-\frac{\tau}{\epsilon}\tilde{f}(\tau)\chi_B(x)\mbox{\boldmath $a$},
$$
where
$$\displaystyle
\tilde{f}(\tau)=\int_0^T e^{-\tau t}f(t)dt.
\tag {1.4}
$$
Note that
$$\displaystyle
\int_0^Te^{-\tau t}\mbox{\boldmath $J$}(x,t)dt
=-\frac{\epsilon}{\tau}\mbox{\boldmath $f$}(x,\tau).
$$

Define
$$
\displaystyle
\mbox{\boldmath $W$}_e(x,\tau)
=\int_0^T e^{-\tau t}\mbox{\boldmath $E$}(x,t)dt,\,\,x\in\Bbb R^3\setminus\overline D.
\tag {1.5}
$$

We call the map defined by
$$\displaystyle
\tau\longmapsto \int_B\mbox{\boldmath $f$}\cdot(\mbox{\boldmath $W$}_e-\mbox{\boldmath $V$})dx
$$
the {\it indicator function}.
The indicator function can be computed from our observation data $\mbox{\boldmath $E$}$ on $B$ over time interval $]0,\,T[$
since we have (1.5).

The following results give us some solutions to the problem raised above.

\proclaim{\noindent Theorem 1.1.}
Assume that $\partial D$ is $C^2$.
Let $f$ satisfy that there exists $\gamma\in\Bbb R$ such that
$$\displaystyle
\liminf_{\tau\longrightarrow\infty}\tau^{\gamma}\vert\tilde{f}(\tau)\vert>0.
\tag {1.6}
$$
If $T>2\sqrt{\mu\epsilon}\text{dist}\,(D,B)$,
then, there exists $\tau_0>0$ such that, for all $\tau\ge\tau_0$
$$\displaystyle
\displaystyle
\int_{B}\mbox{\boldmath $f$}\cdot(\mbox{\boldmath $W$}_e-\mbox{\boldmath $V$})dx>0.
$$
Moreover, we have the following formula:
$$\displaystyle
\lim_{\tau\longrightarrow\infty}\frac{1}{\tau}\log\left\vert\int_{B}\mbox{\boldmath $f$}\cdot
(\mbox{\boldmath $W$}_e-\mbox{\boldmath $V$})dx\right\vert=
-2\sqrt{\mu\epsilon}\text{dist}\,(D,B).
\tag {1.7}
$$

\endproclaim

A remarkable point in this theorem is: there is no restriction on direction $\mbox{\boldmath $a$}$ in (1.2).
Define $d_{\partial D}(p)=\inf_{x\in\partial D}\vert x-p\vert$
and $B_{d_{\partial D}(p)}(p)=\{x\in\Bbb R^3\,\vert\,\vert x-p\vert<d_{\partial D}(p)\}$.
Since we have $\text{dist}\,(D,B)=d_{\partial D}(p)-\eta$,
we can find the sphere $\partial B_{d_{\partial D}(p)}(p)$ via (1.7) regardless of the direction
of $\mbox{\boldmath $a$}$ at any time.  This sphere is the maximum one whose exterior encloses the unknown obstacle.

As is introduced in the author's previous papers \cite{IE3, IE4,
IE5} we denote by $\Lambda_{\partial D}(p)$ the set $\partial
D\cap\partial B_{d_{\partial D}(p)}(p)$.  We call this set the
{\it first reflector} from $p$ to $\partial D$ and the points in
the first reflector are called the {\it first reflection points},
going from $p$ to $\partial D$. Using Theorem 1.1, one can also give a
criterion for a given direction $\omega\in S^2$ whether the point
$p+d_{\partial D}(p)\omega$ belongs to $\partial D$ since as
pointed out in \cite{IE3, IE4, IE5, IH} we have: if $p+d_{\partial
D}(p)\omega$ belongs to $\partial D$, then $d_{\partial
D}(p+sd_{\partial D}(p)\omega)=d_{\partial D}(p)-s$; if
$p+d_{\partial D}(p)\omega$ does not belong to $\partial D$, then
$d_{\partial D}(p+sd_{\partial D}(p)\omega)>d_{\partial D}(p)-s$.
Here $s\in\,]0,\,1[$ and is fixed.  Note that one can always
compute $d_{\partial D}(p+sd_{\partial D}(p)\omega)$ via (1.7)
using a suitable input current supported around $p+sd_{\partial
D}(p)\omega$ and the electronic wave observed at the same place as
the support of the current.

Thus, we obtain the following result which makes use of {\it infinitely} many electromagnetic waves corresponding
to {\it infinitely} many input sources.

\proclaim{\noindent Corollary 1.1.}  Let $p\in\Bbb R^3\setminus\overline D$.
Assume that $d_{\partial D}(p)$ is known.  Fix $\mbox{\boldmath $a$}$,
$\delta\in\,]0,\,d_{\partial D}(p)[$, $\eta'\in\,]0,\,d_{\partial D}(p)-\delta[$
and $f$ satisfying (1.6) for a $\gamma\in\Bbb R$.
Let $T$ satisfy
$$\displaystyle
T>2\sqrt{\mu\epsilon}\,\sup_{p'\in \partial B_{\delta}(p)}\text{dist}\,(D,B_{\eta'}(p')).
$$

Then, one can extract $\Lambda_{\partial D}(p)$ itself from $\mbox{\boldmath $E$}(x,t)$ given at all $x\in B_{\eta'}(p')$,
$t\in\,]0,\,T[$ and $p'\in\partial B_{\delta}(p)$ for $\mbox{\boldmath $J$}$ given by (1.2) where $f$ is as above
and $B$ replaced with $B_{\eta'}(p')$.

\endproclaim

It would be interesting to find a {\it constructive} and {\it
exact} method for extracting $\Lambda_{\partial D}(p)$ itself from
a {\it single} electromagnetic wave, however, at the present time,
we have only a positive result for a scalar wave equation with
Dirichlet boundary condition on the boundary of the obstacle
\cite{IE4}.  The point is to make use of the observed data
restricted to infinitely many closed balls contained in $B$ for a
fixed initial data supported on $\overline B$, that is the
so-called {\it bistatic data}.

The condition (1.6) is a restriction on the strength of the source at $t=0$.
Note that we have $\displaystyle
\tilde{f}(\tau)=O(\tau^{-3/2})$ as $\tau\longrightarrow\infty$ since $f\in H^1(0,\,T)$ and $f(0)=0$.
Thus, $\gamma$ in (1.6) has to satisfy $\gamma\ge 3/2$.  For example, any $f\in H^1(0,\,T)$ such that
$f(t)=t\sin\omega t$ for all $t\in]0,\,\epsilon[$ with $0<\epsilon\le T$ and $\omega>0$,
satisfies (1.6) for $\gamma=3$ since, as $\tau\longrightarrow\infty$
$$\displaystyle
\int_0^{\epsilon}e^{-\tau t}t\sin\,\omega tdt
=\frac{2\tau\omega}{(\tau^2+\omega^2)^2}+O(\tau^{-1}e^{-\epsilon\tau}).
$$

Let $q\in\Lambda_{\partial D}(p)$.
Let $S_q(\partial D)$ and $S_q(\partial
B_{d_{\partial D}(p)}(p))$ denote the {\it shape operators} (or {\it Weingarten maps}) at $q$
of $\partial D$ and $\partial B_{d_{\partial D}(p)}(p)$ with
respect to $\mbox{\boldmath $\nu$}_q$ and $-\mbox{\boldmath $\nu$}_q$, respectively
(see \cite{O} for the notion of the shape operator).
These are symmetric linear operators on the common tangent space $T_q\partial
D=T_q\partial B_{d_{\partial D}(p)}(p)$. We have always
$S_q(\partial B_{d_{\partial D}(p)}(p))-S_q(\partial D)\ge 0$
since $q$ attains the minimum value of the function $\partial D\ni
y\longmapsto \vert y-p\vert$. In general, given $p$ the first
reflector from $p$ to $\partial D$ can be an infinite set, even
more, a continuum. For example, imagine the case when a part of
$\partial D$ coincides with that of $\partial B_{d_{\partial
D}(p)}(p)$. Note also that, in that case, we have $S_q(\partial
B_{d_{\partial D}(p)}(p))=S_q(\partial D)$ the points $q$ in that
part.

\proclaim{\noindent Theorem 1.2.}
Assume that $\partial D$ is $C^4$; $\Lambda_{\partial D}(p)$ is finite and satisfies
$$\displaystyle
\text{det}\,(S_q(\partial B_{d_{\partial D}(p)}(p))-S_q(\partial D))>0,\,\,\forall q\in\Lambda_{\partial D}(p).
\tag {1.8}
$$
Moreover, assume that
$$\displaystyle
\exists q\in\Lambda_{\partial D}(p)\,\,\,\,\vert\mbox{\boldmath $a$}\cdot\mbox{\boldmath $\nu$}_q\vert\not=1.
\tag {1.9}
$$

Let $f$ satisfy (1.6) for a $\gamma\in\Bbb R$.
If $T>2\sqrt{\mu\epsilon}\text{dist}\,(D,B)$, then we have
$$\begin{array}{c}
\displaystyle
\lim_{\tau\longrightarrow\infty}
\frac{\displaystyle
\tau^2e^{2\tau\sqrt{\mu\epsilon}\text{dist}\,(D,B)}}{\displaystyle\tilde{f}(\tau)^2}
\int_{B}\mbox{\boldmath $f$}\cdot(\mbox{\boldmath $W$}_e-\mbox{\boldmath $V$})dx\\
\\
\displaystyle
=\frac{\pi}{2\epsilon^2}
\left(\frac{\eta}{d_{\partial D}(p)}\right)^2\sum_{q\in\Lambda_{\partial D}(p)}
\frac{1-(\mbox{\boldmath $a$}\cdot\mbox{\boldmath $\nu$}_q)^2}
{\sqrt{\text{det}\,(S_q(\partial B_{d_{\partial D}(p)}(p))-S_q(\partial D))}}.
\end{array}
\tag {1.10}
$$

\endproclaim

(1.9) is a restriction on the direction of $\mbox{\boldmath $a$}$ in (1.2).
Since $\vert\mbox{\boldmath $a$}\cdot\mbox{\boldmath $\nu$}_q\vert=1$ if and only if $\mbox{\boldmath $a$}=\pm\,\mbox{\boldmath $\nu$}_q$,
(1.9) means that there is no first reflection point from $p$ on the straight line
passing through $p$ and parallel
to $\mbox{\boldmath $a$}$.
It is clear that if $\Lambda_{\partial D}(p)$ consists of at least {\it three points}, then (1.9) is satisfied.

The denominator $\sqrt{\text{det}\,(S_q(\partial B_{d_{\partial
D}(p)}(p))-S_q(\partial D))}$ in the right-hand side on (1.10) is
independent of $\mbox{\boldmath $a$}$ and numerator
$1-(\mbox{\boldmath $a$}\cdot\mbox{\boldmath $\nu$}_q)^2$ becomes
maximum when $\mbox{\boldmath $a$}$ is perpendicular to
$\mbox{\boldmath $\nu$}_q$; small when $\mbox{\boldmath $a$}\times\mbox{\boldmath $\nu$}_q\approx\mbox{\boldmath $0$}$.
Thus, formula (1.10) shows us an effect of the {\it directivity} of the source term on
extracting information about the geometry of the unknown obstacle from the observation data.
Someone may think that this fact has similarity to a well known fact in the dipole antenna
theory(e.g., \cite{B}), that is, the maximum radiation from the antenna is directed
along right angles to the dipole.

Note also that since we have
$$\displaystyle
\int_B\mbox{\boldmath $f$}\cdot\mbox{\boldmath $W$}_edx
=-\frac{\tau}{\epsilon}\tilde{f}(\tau)\int_B\mbox{\boldmath $a$}\cdot \mbox{\boldmath $W$}_edx,
$$
from (1.5) we know that in (1.7) and (1.10) instead of all the components of $\mbox{\boldmath $E$}$ we need only
$\mbox{\boldmath $a$}\cdot\mbox{\boldmath $E$}$.

It is a due course to deduce the following corollary from Theorem
1.2 (see \cite{IE3, IE5}).

\proclaim{\noindent Corollary 1.2.} Assume that $\partial D$ is
$C^4$.  Let $p\in\Bbb R^3\setminus\overline D$ and assume that
$q\in\Lambda_{\partial D}(p)$ is known. Let $B_1$ and $B_2$ denote
two open balls centred at $p-s_j(p-q)/\vert p-q\vert$, $j=1,2$,
respectively with $0<s_1<s_2<\vert p-q\vert$ and satisfy
$\overline B_1\cup\overline B_2\subset\Bbb R^3\setminus\overline
D$. Let $\mbox{\boldmath $J$}_j$ be the $\mbox{\boldmath $J$}$
given by (1.2) in which $B=B_j$ and $f(t)=f_j(t)$ satisfying (1.6)
for a $\gamma=\gamma_j\in\Bbb R$; $\mbox{\boldmath $E$}_j$ with
$j=1,2$ denote the corresponding electric fields governed by
(1.1).

If $T>2\sqrt{\mu\epsilon}\max_{j=1,2}\text{dist}\,(D, B_j)$ and
$\mbox{\boldmath $a$}\times(p-q)\not=\mbox{\boldmath $0$}$, then
one can extract the Gauss curvature $K_{\partial D}(q)$ of $\partial D$ at $q$ and mean curvature
$H_{\partial D}(q)$ with respect to $\mbox{\boldmath $\nu$}_q$ from
$\mbox{\boldmath $a$}\cdot\mbox{\boldmath $E$}_j$ on $B_j$ with $j=1,2$ over time interval $]0,\,T[$.

\endproclaim

Note that $\mbox{\boldmath $\nu$}_q=(p-q)/\vert p-q\vert$ for $q\in\Lambda_{\partial D}(p)$.
Thus $\mbox{\boldmath $a$}\times(p-q)\not=\mbox{\boldmath $0$}$ if and only if $\vert\mbox{\boldmath $a$}\cdot\mbox{\boldmath $\nu$}_q\vert\not=1$.

Briefly speaking, Corollary 1.2 says that: one can completely know
the Gauss and mean curvatures of the boundary of the obstacle at a known first reflection point $q$,
going from a given point $p$ outside the obstacle by observing
{\it two} reflected electric fields generated by two sources
whose centres are placed on the segment connecting $p$ and $q$.  Thus, one can know an approximate shape
of the boundary of unknown obstacle at a known first reflection point by using two electromagnetic waves.

The concrete procedure in Corollary 1.2 for extracting both the Gauss and mean curvatures at a known first reflection point
consists of the following steps.

(i)  Compute $R_j$ with $j=1,2$ given by
$$\displaystyle
R_j=\lim_{\tau\longrightarrow\infty}
\frac{\displaystyle
\tau^2e^{2\tau\sqrt{\mu\epsilon}(d_{\partial D}(p)-2s_j)}}{\displaystyle\tilde{f_j}(\tau)^2}
\int_{ B_j}
\mbox{\boldmath $f$}_j\cdot(\mbox{\boldmath $W$}_e^j-\mbox{\boldmath $V$}_j)dx,
$$
where $j$ in $\mbox{\boldmath $f$}_j$, $\mbox{\boldmath $W$}_e^j$ and $\mbox{\boldmath $V$}_j$
indicates that they are the $\mbox{\boldmath $f$}$, $\mbox{\boldmath $W$}_e$ and $\mbox{\boldmath $V$}$
corresponding to $f_j$, $\mbox{\boldmath $E$}_j$ and $\mbox{\boldmath $J$}_j$ in a trivial manner.

(ii)  Compute $X_j$ with $j=1,2$ given by
$$\displaystyle
X_j=\left(\frac{1-(\mbox{\boldmath $a$}\cdot\mbox{\boldmath $\nu$}_q)^2}{R_j}\right)^2
\left\{\frac{\pi}{2\epsilon^2}\left(\frac{\eta_j}{d_{\partial D}(p)-s_j}\right)^2\right\}^2.
$$

(iii) Solve the following linear equations with unknowns $Y_1$ and $Y_2$:
$$\displaystyle
\left(\begin{array}{lr}
\displaystyle
-2\lambda_1  & 1\\
\\
\displaystyle
-2\lambda_2 & 1
\end{array}
\right)
\left(
\begin{array}{c}
\displaystyle
Y_1\\
\\
\displaystyle
Y_2
\end{array}
\right)
=\left(
\begin{array}{c}
\displaystyle
X_1\\
\\
\displaystyle
X_2
\end{array}
\right)
-
\left(
\begin{array}{c}
\displaystyle
\lambda_1^2\\
\\
\displaystyle
\lambda_2^2
\end{array}
\right),
$$
where $\lambda_j=(d_{\partial D}(p)-s_j)^{-1}$ with $j=1,2$.

Then, we obtain the Gauss and mean curvatures at $q$ by the formulae: $K_{\partial D}(p)=Y_2$ and $H_{\partial D}(q)=Y_1$.
Note that we have made use of the following trivial facts and formulae as pointed out
in \cite{IE3, IE5}:

$\bullet$  $\Lambda_{\partial D}(p-s_j\mbox{\boldmath $\nu$}_q)=\{q\}$ for $q\in\Lambda_{\partial D}(p)$
and $S_q(\partial B_{d_{\partial D}(p-s_j\mbox{\boldmath $\nu$}_q)}
(p-s_j\mbox{\boldmath $\nu$}_q))-S_q(\partial B_{d_{\partial D}(p)}(p))$ is {\it positive definite} on the common tangent space
at $q$;

$\bullet$  $\text{dist}\,(D,B_j)=d_{\partial D}(p)-2s_j$ and $d_{\partial D}(p-s_j\mbox{\boldmath $\nu$}_q)
=d_{\partial D}(p)-s_j$ for $q\in\Lambda_{\partial D}(p)$;

$\bullet$  $\text{det}\,(\lambda I-S_q(\partial D))=\lambda^2-2\lambda H_{\partial D}(q)+K_{\partial D}(q)$;

$\bullet$  $S_q(\partial B_{d_{\partial D}(p-s_j\mbox{\boldmath $\nu$}_q)}(p-s_j\mbox{\boldmath $\nu$}_q))=\lambda_j I$.

We think that Corollary 1.2 shows us an advantage of the {\it near
field} measurement.  For this, note that in the third step above $\lambda_j\longrightarrow 0$ as $d_{\partial D}(p)\longrightarrow\infty$
and thus one cannot find $Y_1$.
Compare also the results with those of \cite{MT}
where the information about the mean curvature never appear
explicitly in the {\it scattering kernel} which is the observation
data in the context of the Lax-Phillips scattering theory
\cite{LP}.

Note that our result can be applied to a cavity inside a large obstacle
which is connected with its exterior by a borehole.
This is the case when $D$ {\it encloses} almost $B$.
In this case
it is not suitable to use an {\it infinitely extended plane wave} as an approximation of the incident wave unlike \cite{MT}.
See also \cite{AMS} for some comments on the comparison between incident plane and spherical waves in the frequency domain.

The outline of this paper is as follows.
To study the asymptotic behaviour of the indicator function as
$\tau\longrightarrow\infty$ we need some preliminary facts about
$\mbox{\boldmath $V$}$. In Section 2, using the mean value theorem
for the modified Helmholtz equation, we give an explicit
computation formula for $\mbox{\boldmath $V$}$ outside $B$. This
formula is found in Subsection 2.1 and enables us to study the
asymptotic behaviour of the energy integral $J(\tau)$ of $\mbox{\boldmath
$V$}$ over $D$ as $\tau\longrightarrow\infty$ in Subsection 2.2,
where
$$\displaystyle
J(\tau)
=\frac{1}{\epsilon\mu}
\int_D\vert\nabla\times\mbox{\boldmath $V$}\vert^2dx
+\tau^2\int_{D}\vert\mbox{\boldmath $V$}\vert^2dx.
\tag {1.11}
$$
However, unlike the previous applications of the enclosure method to scalar wave
equations (see e.g., \cite{IE5}), we need an {\it upper bound} of $L^2$-norm of the
Jacobian matrix $\mbox{\boldmath $V$}'$ over $D$ in terms of $J(\tau)$.
This is not trivial and described in Subsection 2.2.

Theorem 1.1 is proved in Section 3.  The proof is based on a {\it brief} asymptotic formula
of the indicator function and the resulted upper and lower bound in terms of $J(\tau)$.

Theorem 1.2 is proved in Section 4. The proof is
based on the {\it precise} asymptotic formula of the indicator
function stated in Theorem 4.1 and the leading profile of $J(\tau)$ obtained in Section 2 via the Laplace method.
The precise asymptotic formula is derived from a combination of the brief
asymptotic formula of the indicator function and the {\it asymptotic coincidence} of $J(\tau)$ with $E(\tau)$ defined by
$$\displaystyle
E(\tau)
=\frac{1}{\epsilon\mu}
\int_{\Bbb R^3\setminus\overline D}\vert\nabla\times(\mbox{\boldmath $W$}_e-\mbox{\boldmath $V$})\vert^2dx
+\tau^2\int_{\Bbb R^3\setminus\overline D}\vert\mbox{\boldmath $W$}_e-\mbox{\boldmath $V$}\vert^2dx.
\tag {1.12}
$$
The proof of the asymptotic coincidence of $E(\tau)$ and $J(\tau)$ is based on the
{\it reflection principle} across curved surface $\partial D$ for the Maxwell system as
described in Propositions 4.1 and 4.2 and a representation of
$E(\tau)-J(\tau)$ in terms of the reflection which are trivial for a scalar wave equation case.
Then, we apply the Lax-Phillips
reflection argument \cite{LP} to the difference. This story is
parallel to the previous scalar wave equation cases \cite{IE3,
IE4, IE5}, however, a proper problem for system of partial
differential equations occurs in proving Theorem 4.1.
In order to apply their argument, we need an upper bound of the $L^2$-norm of
the Jacobian matrix of the so-called reflected solution $\mbox{\boldmath $W$}_e-\mbox{\boldmath $V$}$
in terms of $E(\tau)$.  However, it seems
difficult to obtain such an estimate directly and instead, we give the
upper bound in terms of $J(\tau)$ directly.
This way is different from the original
Lax-Phillips reflection argument and makes the argument for the
proof of the asymptotic coincidence of $E(\tau)$ with $J(\tau)$
straightforward compared with the scalar wave equation case.

In Appendix we describe some differential identities for the
vector fields obtained by the reflection across $\partial D$ and
the resulted reflection formula described in Proposition 4.2 is
proved. Note that the regularity assumption that $\partial D$ is
$C^4$ in Theorem 1.2 is more restrictive compared with the scalar
wave equation case \cite{IE3, IE4, IE5} in which the corresponding
theorems are valid for $C^3$-smooth boundary.  This is coming from
the difference of the reflection principle used.  Therein only a
change of {\it independent} variables is used, however, for
Maxwell's equations, the reflection principle involves also a
change of {\it dependent} variables and this requires a higher
regularity.

\section{Preliminary facts about $\mbox{\boldmath $V$}$}

In this section first we give a detailed expression of
$\mbox{\boldmath $V$}$. Then using the expression we give an
asymptotic behaviour of some integrals involving $\mbox{\boldmath
$V$}$.

\subsection{An explicit form of $\mbox{\boldmath $V$}$ outside of $B$}

Here we give an explicit computation formula of the weak solution of (1.3) in $\Bbb R^3\setminus B$.
First, assume that $\mbox{\boldmath $V$}$ has the form
$$\displaystyle
\mbox{\boldmath $V$}=\mbox{\boldmath $V$}_0+\mbox{\boldmath $V$}_1,
$$
where $\mbox{\boldmath $V$}_0$ and $\mbox{\boldmath $V$}_1$ are two vector-valued functions on the whole space.

Write
$$\begin{array}{c}
\displaystyle
\frac{1}{\mu\epsilon}\nabla\times\nabla\times\mbox{\boldmath $V$}+\tau^2\mbox{\boldmath $V$}+\mbox{\boldmath $f$}(x,\tau)\\
\\
\displaystyle
=\left\{-\frac{1}{\mu\epsilon}(\triangle-\mu\epsilon\tau^2)\mbox{\boldmath $V$}_0+\mbox{\boldmath $f$}(x,\tau)\right\}
+\left(\tau^2\mbox{\boldmath $V$}_1+\frac{1}{\mu\epsilon}\nabla(\nabla\cdot\mbox{\boldmath $V$}_0)\right)
+\frac{1}{\mu\epsilon}\nabla\times\nabla\times\mbox{\boldmath $V$}_1.
\end{array}
$$
From this we see that if
$$\displaystyle
-\frac{1}{\mu\epsilon}(\triangle-\mu\epsilon\tau^2)\mbox{\boldmath $V$}_0+\mbox{\boldmath $f$}(x,\tau)=\mbox{\boldmath $0$}
\tag {2.1}
$$
and
$$\displaystyle
\tau^2\mbox{\boldmath $V$}_1+\frac{1}{\mu\epsilon}\nabla(\nabla\cdot\mbox{\boldmath $V$}_0)=\mbox{\boldmath $0$},
$$
then $\nabla\times\mbox{\boldmath $V$}_1=\mbox{\boldmath $0$}$ and thus $\mbox{\boldmath $V$}=\mbox{\boldmath $V$}_0+\mbox{\boldmath $V$}_1$
satisfies (1.3) formally.

From this formal argument we have the following construction of the weak solution of
(1.3) for general $\mbox{\boldmath $f$}(\,\cdot\,,\tau)\in L^2(\Bbb R^3)^3$ such that
$\displaystyle\text{supp}\,\mbox{\boldmath $f$}(\,\cdot\,,\tau)\subset\overline B$.

Let $\mbox{\boldmath $V$}_0=\mbox{\boldmath $V$}_0(\,\cdot\,,\tau)\in\,H^1(\Bbb R^3)^3$ be the unique weak solution of (2.1).
It is well known that $\mbox{\boldmath $V$}_0$ has the form
$$\displaystyle
\mbox{\boldmath $V$}_0=\mbox{\boldmath $V$}_0(x,\tau)
=-\frac{\mu\epsilon}{4\pi}
\int_{B}\frac{e^{-\sqrt{\mu\epsilon}\tau\vert x-y\vert}}{\vert x-y\vert}\mbox{\boldmath $f$}(y,\tau)dy.
\tag {2.2}
$$
Then, for each fixed $\tau$ by the interior regularity or from the expression
we see that $\mbox{\boldmath $V$}_0\in H^2_{\text{loc}}\,(\Bbb R^3)^3$; $\mbox{\boldmath $V$}_0$ is smooth
outside $B$; $\mbox{\boldmath $V$}_0$ together with its all derivatives are exponentially decaying as $\vert x\vert\longrightarrow\infty$.
Thus we have $\mbox{\boldmath $V$}_0\in H^2(\Bbb R^3)^3$.

Define $\mbox{\boldmath $V$}_1=\mbox{\boldmath $V$}_1(\,\cdot\,,\tau)\in L^2(\Bbb R^3)^3$ by the formula
$$\displaystyle
\mbox{\boldmath $V$}_1=-\frac{1}{\tau^2\mu\epsilon}\nabla(\nabla\cdot\mbox{\boldmath $V$}_0).
\tag {2.3}
$$
$\mbox{\boldmath $V$}_1$ is also smooth outside $B$ and, for each fixed $\tau$, $\mbox{\boldmath $V$}_1$ together with all the derivatives
are exponentially decaying as $\vert x\vert\longrightarrow\infty$; $\nabla\times\mbox{\boldmath $V$}_1=\mbox{\boldmath $0$}$ in $\Bbb R^3$.

Then, $\mbox{\boldmath $V$}=\mbox{\boldmath $V$}_0+\mbox{\boldmath $V$}_1\in L^2(\Bbb R^3)^3$ satisfies
$\nabla\times\mbox{\boldmath $V$}=\nabla\times\mbox{\boldmath $V$}_0\in L^2(\Bbb R^3)^3$.  Thus
we have $\mbox{\boldmath $V$}\in H(\text{curl}\,,\Bbb R^3)$.
It is easy to see that this $\mbox{\boldmath $V$}$ satisfies (1.3) in the weak sense.  Thus, by the uniqueness
of the weak solution of (1.3) we conclude that the weak solution of (1.3) has the expression
$$\displaystyle
\mbox{\boldmath $V$}=\mbox{\boldmath $V$}_0+\mbox{\boldmath $V$}_1,
\tag {2.4}
$$
where $\mbox{\boldmath $V$}_0$ and $\mbox{\boldmath $V$}_1$ are given by (2.2) and (2.3), respectively.
Note that this argument for the construction of $\mbox{\boldmath $V$}$
is based on the form of the fundamental solution for the operator
$\displaystyle (1/\mu\epsilon)\nabla\times\nabla\times\,\cdot\,-k^2\,\cdot\,$ with $k>0$ (e.g., see \cite{ABF}).

Let $x\in\Bbb R^3\setminus\overline B$.
In what follows, we omit to indicate the dependence of several functions of $x$ on the parameter $\tau$.
By the mean value theorem for the modified Helmholtz equation \cite{CH}, we know that
$$\displaystyle
\frac{1}{4\pi}\int_B\frac{e^{-\tau\sqrt{\mu\epsilon}\vert x-y\vert}}{\vert x-y\vert}dy
=\frac{\varphi(\tau\sqrt{\mu\epsilon}\eta)}{(\tau\sqrt{\mu\epsilon})^3}\frac{e^{-\tau\sqrt{\mu\epsilon}\vert x-p\vert}}{\vert x-p\vert},
$$
where
$\displaystyle
\varphi(\xi)=\xi\cosh\,\xi-\sinh\,\xi$.
Thus $\mbox{\boldmath $V$}_0$ given by (2.2) takes the form
$$\displaystyle
\mbox{\boldmath $V$}_0(x)
=K(\tau)\tilde{f}(\tau)v(x)\mbox{\boldmath $a$},
\tag {2.5}
$$
where $\tilde{f}(\tau)$ is given by (1.4),
$$\displaystyle
v(x)=\frac{e^{-\tau\sqrt{\mu\epsilon}\vert x-p\vert}}{\vert x-p\vert}
$$
and
$$\displaystyle
K(\tau)=\frac{\mu\tau\varphi(\tau\sqrt{\mu\epsilon}\eta)}{(\tau\sqrt{\mu\epsilon})^3}.
$$

A straightforward computation gives
$$\displaystyle
\nabla v(x)
=-\left(\frac{\tilde{\tau}}{\vert x-p\vert}
+\frac{1}{\vert x-p\vert^2}\right)(x-p)v(x)
$$
and
$$\begin{array}{c}
\displaystyle
(\nabla v)'(x)\\
\\
\displaystyle
=v(x)\left\{
\left(\frac{\tilde{\tau}^2}{\vert x-p\vert^2}
+\frac{3\tilde{\tau}}{\vert x-p\vert^3}
+\frac{3}{\vert x-p\vert^4}\right)(x-p)\otimes(x-p)
-\left(\frac{\tilde{\tau}}{\vert x-p\vert}+\frac{1}{\vert x-p\vert^2}\right)I_3
\right\},
\end{array}
$$
where $\tilde{\tau}=\tau\sqrt{\mu\epsilon}$.
Since
$\displaystyle
\nabla(\nabla\cdot\mbox{\boldmath $V$}_0)
=K(\tau)\tilde{f}(\tau)(\nabla v)'\mbox{\boldmath $a$}$,
it follows from (2.3) that
$$\begin{array}{c}
\displaystyle
\mbox{\boldmath $V$}_1(x)
=-\frac{K(\tau)\tilde{f}(\tau)}{\tilde{\tau}^2}v(x)\times\\
\\
\displaystyle
\left\{
\left(\tilde{\tau}^2
+\frac{3\tilde{\tau}}{\vert x-p\vert}
+\frac{3}{\vert x-p\vert^2}\right)\frac{x-p}{\vert x-p\vert}\otimes\frac{x-p}{\vert x-p\vert}
-\left(\frac{\tilde{\tau}}{\vert x-p\vert}+\frac{1}{\vert x-p\vert^2}\right)I_3\right\}\mbox{\boldmath $a$}.
\end{array}
$$
Substituting this and (2.5) into (2.4), we obtain the following explicit formula
of $\mbox{\boldmath $V$}$ outside $B$:
$$\displaystyle
\mbox{\boldmath $V$}(x)
=K(\tau)\tilde{f}(\tau)v(x)\mbox{\boldmath $M$}(x;p)\mbox{\boldmath $a$},
\tag {2.6}
$$
where
$$\displaystyle
\mbox{\boldmath $M$}(x;p)
=AI_3-B\frac{x-p}{\vert x-p\vert}\otimes\frac{x-p}{\vert x-p\vert},
\tag {2.7}
$$
$$\displaystyle
A=A(x,\tau)=1+\frac{1}{\tau\sqrt{\mu\epsilon}}
\left(\frac{1}{\vert x-p\vert}
+\frac{1}{\tau\sqrt{\mu\epsilon}\vert x-p\vert^2}\right)
$$
and
$$\displaystyle
B=B(x,\tau)=1+\frac{3}{\tau\sqrt{\mu\epsilon}}
\left(\frac{1}{\vert x-p\vert}
+\frac{1}{\tau\sqrt{\mu\epsilon}\vert x-p\vert^2}\right).
$$

\subsection{Two basic lemmas about $J(\tau)$ and $\mbox{\boldmath $V$}'$}

Let $\overline B\subset\Bbb R^3\setminus\overline D$.  The following two lemmas are concerned with the asymptotic behaviour
of some integrals involving $\mbox{\boldmath $V$}$ and its derivatives over $D$ which is one of the key points in this paper.

\proclaim{\noindent Lemma 2.1.}

(i)  We have
$$\displaystyle
\limsup_{\tau\longrightarrow\infty}\tau^{3}e^{\displaystyle 2\tau\sqrt{\mu\epsilon}\text{dist}\,(D,B)}J(\tau)<\infty.
\tag {2.8}
$$

(ii) Assume that $\partial D$ is Lipschitz.  Let $f$ satisfy (1.6) for a $\gamma\in\Bbb R$.
We have
$$\displaystyle
\liminf_{\tau\longrightarrow\infty}\tau^{5+2\gamma}e^{\displaystyle 2\tau\sqrt{\mu\epsilon}\text{dist}\,(D,B)}J(\tau)>0.
\tag {2.9}
$$
\endproclaim

{\it\noindent Proof.}
For convenience we introduce
$\displaystyle
\tilde{\tau}=\tau\sqrt{\mu\epsilon}$.
Since $\vert\tilde{f}(\tau)\vert=O(\tau^{-3/2})$
and
$$\displaystyle
\varphi(\tilde{\tau}\eta)
=\frac{\tilde{\tau}\eta e^{\tilde{\tau}\eta}}{2}(1+O(\tau^{-1})),
\tag {2.10}
$$
we have $K(\tau)\tilde{f}(\tau)=O(\tau^{-5/2}e^{\tilde{\tau}\eta})$.
Thus, it follows from (2.6) and (2.7) that, for all $x\in D$
$$\displaystyle
\vert\mbox{\boldmath $V$}(x)\vert
\le C\tau^{-5/2}e^{\displaystyle\tilde{\tau}\eta}v(x).
$$
Moreover, since $\nabla\times\mbox{\boldmath $V$}_1=0$, from (2.5)
we have
$$\begin{array}{c}
\displaystyle
\nabla\times\mbox{\boldmath $V$}
=-\tilde{\tau}K(\tau)\tilde{f}(\tau)v(x)
\left(1+\frac{1}{\tilde{\tau}\vert x-p\vert}\right)\frac{x-p}{\vert x-p\vert}\times\mbox{\boldmath $a$}
\end{array}
\tag {2.11}
$$
and thus, for all $x\in D$
$$\displaystyle
\vert\nabla\times\mbox{\boldmath $V$}(x)\vert
\le C\tau^{-3/2}e^{\displaystyle\tilde{\tau}\eta}v(x).
$$
Therefore one gets
$$\displaystyle
\tau^3J(\tau)\le Ce^{\displaystyle 2\tilde{\tau}\eta}\int_D\vert v(x)\vert^2dx.
$$
Now (2.8) is clear since we have
$$\displaystyle
\int_D\vert v(x)\vert^2dx
\le \frac{1}{d_{\partial D}(p)^2}\int_De^{-2\tilde{\tau}\vert x-p\vert}dx
=O(e^{\displaystyle -2\tilde{\tau}d_{\partial D}(p)})
$$
and $d_{\partial D}(p)-\eta=\text{dist}\,(D,B)$.

Next we give a proof of (2.9).
From (2.7) we have
$$\begin{array}{c}
\displaystyle
\vert\mbox{\boldmath $M$}(x;p)\mbox{\boldmath $a$}\vert^2
=\mbox{\boldmath $M$}(x;p)^2\mbox{\boldmath $a$}\cdot\mbox{\boldmath $a$}
\\
\\
\displaystyle
=A^2\vert\mbox{\boldmath $a$}\vert^2+(B^2-2AB)\left\vert\mbox{\boldmath $a$}\cdot\frac{x-p}{\vert x-p\vert}\right\vert^2\\
\\
\displaystyle
=\left(1-\left\vert\mbox{\boldmath $a$}\cdot\frac{x-p}{\vert x-p\vert}\right\vert^2\right)A^2
+(A^2+B^2-2AB)\left\vert\mbox{\boldmath $a$}\cdot\frac{x-p}{\vert x-p\vert}\right\vert^2\\
\\
\displaystyle
=A^2\left\vert\mbox{\boldmath $a$}\times\frac{x-p}{\vert x-p\vert}\right\vert^2
+(A-B)^2\left\vert\mbox{\boldmath $a$}\cdot\frac{x-p}{\vert x-p\vert}\right\vert^2\\
\\
\displaystyle
\ge\frac{1}{\tau^2\mu\epsilon\vert x-p\vert^2}\left\vert\mbox{\boldmath $a$}\times\frac{x-p}{\vert x-p\vert}\right\vert^2
+\frac{4}{\tau^2\mu\epsilon\vert x-p\vert^2}\left\vert\mbox{\boldmath $a$}\cdot\frac{x-p}{\vert x-p\vert}\right\vert^2
\ge\frac{1}{\tau^2\mu\epsilon\vert x-p\vert^2}.
\end{array}
$$
Thus (2.6) gives
$$\displaystyle
\vert\mbox{\boldmath $V$}(x)\vert^2
\ge \tau^{-2}K(\tau)^2\tilde{f}(\tau)^2\frac{v(x)^2}{\mu\epsilon\vert x-p\vert^2}.
\tag {2.12}
$$
By (2.10) we have $K(\tau)\sim \tau^{-1}e^{\tilde{\tau}\eta}\eta/(2\epsilon)$ as $\tau\longrightarrow\infty$.
Then, it follows from (2.12) that
there exist positive constants $C''$ and $\tau_0$ such that, for all $x\in D$ and
$\tau\ge\tau_0$
$$\displaystyle
\vert\mbox{\boldmath $V$}(x)\vert^2
\ge C''\tau^{-4}\tilde{f}(\tau)^2e^{2\tilde{\tau}\eta}v(x)^2
$$
and thus
$$\displaystyle
J(\tau)
\ge C''\tau^{-2}\tilde{f}(\tau)^2e^{\displaystyle 2\tilde{\tau}\eta}\int_{D}\vert v(x)\vert^2 dx.
\tag {2.13}
$$
A standard technique \cite{IH} yields
$$\displaystyle
\liminf_{\tau\longrightarrow\infty}\tau^3e^{\displaystyle 2\tilde{\tau} d_{\partial D}(p)}
\int_{D}\vert v(x)\vert^2 dx>0.
$$
Thus rewriting (2.13) as
$$\displaystyle
e^{\displaystyle 2\tilde{\tau}\text{dist}\,(D,B)}
\tau^{3+2+2\gamma}J(\tau)
\ge C\tau^{2\gamma}\tilde{f}(\tau)^2\times
\tau^3e^{\displaystyle2\tilde{\tau} d_{\partial D}(p)}
\int_{D}\vert v(x)\vert^2 dx,
$$
we obtain (2.9).

\noindent
$\Box$

For the proof of Theorem 1.2 we need a more accurate information
about the asymptotic behaviour of $J(\tau)$ as
$\tau\longrightarrow\infty$.

\proclaim{\noindent Lemma 2.2.}
Let $f$ satisfy: there exists a positive constant $\tau_0$ such that, for all $\tau\ge\tau_0$ $\tilde{f}(\tau)\not=0$.
Assume that $\Lambda_{\partial D}(p)$ is finite and satisfies (1.8).
Then, we have
$$\begin{array}{c}
\displaystyle
\lim_{\tau\longrightarrow\infty}\frac{\tau^2e^{2\tau\sqrt{\mu\epsilon}\text{dist}\,(D,B)}}
{\displaystyle \tilde{f}(\tau)^2}
J(\tau)\\
\\
\displaystyle
=\frac{\pi}{4\epsilon^2}\left(\frac{\eta}{d_{\partial D}(p)}\right)^2
\sum_{q\in\Lambda_{\partial D}(p)}
\frac{1-(\mbox{\boldmath $a$}\cdot\mbox{\boldmath $\nu$}_q)^2}
{\sqrt{\text{det}\,(S_q(\partial B_{d_{\partial D}(p)}(p))-S_q(\partial D))}}.
\end{array}
\tag {2.14}
$$
Moreover, if (1.9) is also satisfied, then, as
$\tau\longrightarrow\infty$
$$\displaystyle
\int_D\vert\mbox{\boldmath $V$}'\vert^2dx=O(J(\tau)).
\tag {2.15}
$$

\endproclaim

{\it\noindent Proof.}
Using (1.3), (1.11) and Integration by parts, we obtain
$$\begin{array}{c}
\displaystyle
J(\tau)=\frac{1}{\mu\epsilon}\int_{\partial D}(\mbox{\boldmath $\nu$}\times\mbox{\boldmath $V$})
\cdot\nabla\times\mbox{\boldmath $V$}dS.
\end{array}
\tag {2.16}
$$
Then, the identity $\displaystyle
(\mbox{\boldmath $\nu$}\times\mbox{\boldmath $V$})\cdot\nabla\times\mbox{\boldmath $V$}
=-\mbox{\boldmath $\nu$}\cdot(\nabla\times\mbox{\boldmath $V$})\times\mbox{\boldmath $V$}
\,\,\text{on}\, \partial D$
yields another expression
$$\displaystyle
J(\tau)=-\int_{\partial D}\mbox{\boldmath $\nu$}\cdot(\nabla\times\mbox{\boldmath $V$})\times\mbox{\boldmath $V$}dS.
\tag {2.17}
$$
Let $x\in\partial D$ and $\displaystyle\mbox{\boldmath $\omega$}_x=(x-p)/\vert x-p\vert$.

We have
$$\begin{array}{c}
\displaystyle
(\mbox{\boldmath $\omega$}_x\times\mbox{\boldmath $a$})
\times(\mbox{\boldmath $a$}-(\mbox{\boldmath $\omega$}_x\cdot\mbox{\boldmath $a$})\mbox{\boldmath $\omega$}_x)\\
\\
\displaystyle
=\mbox{\boldmath $a$}(\mbox{\boldmath $\omega$}_x\cdot\mbox{\boldmath $a$})
-\mbox{\boldmath $\omega$}_x(\mbox{\boldmath $a$}\cdot\mbox{\boldmath $a$})
+(\mbox{\boldmath $\omega$}_x\cdot\mbox{\boldmath $a$})^2\mbox{\boldmath $\omega$}_x
-(\mbox{\boldmath $\omega$}_x\cdot\mbox{\boldmath $a$})\mbox{\boldmath $a$}\\
\\
\displaystyle
=-\mbox{\boldmath $\omega$}_x(\mbox{\boldmath $a$}\cdot\mbox{\boldmath $a$})
+(\mbox{\boldmath $\omega$}_x\cdot\mbox{\boldmath $a$})^2\mbox{\boldmath $\omega$}_x.
\end{array}
$$
Thus one can write
$$\displaystyle
(\mbox{\boldmath $\omega$}_x\times\mbox{\boldmath $a$})
\times(\mbox{\boldmath $a$}-(\mbox{\boldmath $\omega$}_x\cdot\mbox{\boldmath $a$})\mbox{\boldmath $\omega$}_x)
\cdot\mbox{\boldmath $\nu$}_x
=\left(\mbox{\boldmath $m$}(x;p)\mbox{\boldmath $\nu$}_x\right)\mbox{\boldmath $a$}\cdot\mbox{\boldmath $a$},
\tag {2.18}
$$
where
$$\displaystyle
\mbox{\boldmath $m$}(x;p)\mbox{\boldmath $\nu$}_x
=(\mbox{\boldmath $\omega$}_x\cdot\mbox{\boldmath $\nu$}_x)
(\mbox{\boldmath $\omega$}_x
\otimes\mbox{\boldmath $\omega$}_x-I_3).
$$
Let $q\in\Lambda_{\partial D}(p)$. We have $\mbox{\boldmath $\nu$}_q=-\mbox{\boldmath $\omega$}_q$.
Then $\displaystyle
\mbox{\boldmath $m$}(q;p)\mbox{\boldmath $\nu$}_q
=I_3-\mbox{\boldmath $\nu$}_q\otimes\mbox{\boldmath $\nu$}_q$
and hence
$$\displaystyle
(\mbox{\boldmath $m$}(q;p)\mbox{\boldmath $\nu$}_q)\mbox{\boldmath $a$}\cdot\mbox{\boldmath $a$}
=1-(\mbox{\boldmath $a$}\cdot\mbox{\boldmath $\nu$}_q)^2.
\tag {2.19}
$$
Therefore, $\vert\mbox{\boldmath $a$}\cdot\mbox{\boldmath $\nu$}_q\vert\not=1$
if and only if $\displaystyle
(\mbox{\boldmath $m$}(q;p)\mbox{\boldmath $\nu$}_q)\mbox{\boldmath $a$}\cdot\mbox{\boldmath $a$}\not=\mbox{\boldmath $0$}$.

From (2.6) and (2.11) we have
$$\begin{array}{c}
\displaystyle
(\nabla\times\mbox{\boldmath $V$})\times\mbox{\boldmath $V$}\\
\\
\displaystyle
=-\tilde{\tau}K(\tau)^2\tilde{f}(\tau)^2
v(x)^2
\left(1+\frac{1}{\tilde{\tau}\vert x-p\vert}\right)
\left(\mbox{\boldmath $\omega$}_x\times\mbox{\boldmath $a$}\right)\times(\mbox{\boldmath $M$}(x;p)\mbox{\boldmath $a$}),
\end{array}
\tag {2.20}
$$
where $\tilde{\tau}=\tau\sqrt{\mu\epsilon}$.

From (2.7) we have, as $\tau\longrightarrow\infty$ in the compact uniform topology in $\Bbb R^3\setminus\overline B$
$$\displaystyle
\mbox{\boldmath $M$}(x;p)
=I_3-\mbox{\boldmath $\omega$}_x\otimes\mbox{\boldmath $\omega$}_x
+O\left(\frac{1}{\tau}\right)
\tag {2.21}
$$
and this yields
$$\begin{array}{c}
\displaystyle
\left(1+\frac{1}{\tilde{\tau}\vert x-p\vert}\right)
\left(\mbox{\boldmath $\omega$}_x\times\mbox{\boldmath $a$}\right)\times(\mbox{\boldmath $M$}(x;p)\mbox{\boldmath $a$})\\
\\
\displaystyle
=\left(\mbox{\boldmath $\omega$}_x\times\mbox{\boldmath $a$}\right)\times
\left\{\left(I_3-\mbox{\boldmath $\omega$}_x\otimes\mbox{\boldmath $\omega$}_x\right)\mbox{\boldmath $a$}\right\}
+O\left(\frac{1}{\tau}\right)
\end{array}
$$
uniformly for $x\in\partial D$.
Applying this together with (2.18) to (2.20),
we obtain, as $\tau\longrightarrow\infty$
$$\begin{array}{c}
\displaystyle
-\mbox{\boldmath $\nu$}\cdot(\nabla\times\mbox{\boldmath $V$})\times\mbox{\boldmath $V$}
=\tilde{\tau}K(\tau)^2\tilde{f}(\tau)^2
v(x)^2
\left\{(\mbox{\boldmath $m$}(x;p)\mbox{\boldmath $\nu$}_x)\mbox{\boldmath $a$}\cdot\mbox{\boldmath $a$}
+O\left(\frac{1}{\tau}\right)\right\}
\end{array}
$$
uniformly for $x\in\partial D$.
Thus, we have
$$\begin{array}{c}
\displaystyle
-\int_{\partial D}\mbox{\boldmath $\nu$}\cdot(\nabla\times\mbox{\boldmath $V$})\times\mbox{\boldmath $V$}dS\\
\\
\displaystyle
=\tilde{\tau}K(\tau)^2\tilde{f}(\tau)^2
\left\{\int_{\partial D}v(x)^2(\mbox{\boldmath $m$}(x;p)\mbox{\boldmath $\nu$}_x)\mbox{\boldmath $a$}\cdot\mbox{\boldmath $a$}dS
+O\left(\frac{1}{\tau}\right)
\int_{\partial D}v(x)^2dS\right\}.
\end{array}
\tag {2.22}
$$
Under the finiteness of $\Lambda_{\partial D}(p)$ and (1.8), using the Laplace method \cite{BH} we obtain
$$\begin{array}{c}
\displaystyle
\lim_{\tau\longrightarrow\infty}
\tilde{\tau}e^{2\tilde{\tau}d_{\partial D}(p)}\int_{\partial D}\frac{e^{-2\tilde{\tau}\vert x-p\vert}}{\vert x-p\vert^2}(\mbox{\boldmath $m$}(x;p)\mbox{\boldmath $\nu$}_x)\mbox{\boldmath $a$}\cdot\mbox{\boldmath $a$}dS\\
\\
\displaystyle
=\frac{\pi}{d_{\partial D}(p)^2}\sum_{q\in\Lambda_{\partial D}(p)}
\frac{(\mbox{\boldmath $m$}(q;p)\mbox{\boldmath $\nu_q$})\mbox{\boldmath $a$}\cdot\mbox{\boldmath $a$}}
{\sqrt{\text{det}\,S_q(\partial B_{d_{\partial D}(p)}(p))-S_q(\partial D))}}
\end{array}
\tag {2.23}
$$
and
$$
\begin{array}{c}
\displaystyle
\lim_{\tau\longrightarrow\infty}
\tilde{\tau}e^{2\tilde{\tau}d_{\partial D}(p)}\int_{\partial D}\frac{e^{-2\tilde{\tau}\vert x-p\vert}}{\vert x-p\vert^2}dS
=\frac{\pi}{d_{\partial D}(p)^2}\sum_{q\in\Lambda_{\partial D}(p)}
\frac{1}
{\sqrt{\text{det}\,S_q(\partial B_{d_{\partial D}(p)}(p))-S_q(\partial D))}}.
\end{array}
\tag {2.24}
$$
Thus, applying (2.23) and (2.24) to (2.22), we obtain
$$\begin{array}{c}
\displaystyle
-\lim_{\tau\longrightarrow\infty}\frac{e^{2\tilde{\tau}d_{\partial D}(p)}}
{K(\tau)^2\tilde{f}(\tau)^2}
\int_{\partial D}\mbox{\boldmath $\nu$}\cdot(\nabla\times\mbox{\boldmath $V$})\times\mbox{\boldmath $V$}dS\\
\\
\displaystyle
=\frac{\pi}{d_{\partial D}(p)^2}\sum_{q\in\Lambda_{\partial D}(p)}
\frac{(\mbox{\boldmath $m$}(q;p)\mbox{\boldmath $\nu_q$})\mbox{\boldmath $a$}\cdot\mbox{\boldmath $a$}}
{\sqrt{\text{det}\,S_q(\partial B_{d_{\partial D}(p)}(p))-S_q(\partial D))}}.
\end{array}
$$
Thus, from this, (2.17) and (2.19) we obtain
$$\begin{array}{c}
\displaystyle
\lim_{\tau\longrightarrow\infty}\frac{e^{2\tilde{\tau}d_{\partial D}(p)}}
{K(\tau)^2\tilde{f}(\tau)^2}
J(\tau)
=\frac{\pi}{d_{\partial D}(p)^2}
\sum_{q\in\Lambda_{\partial D}(p)}
\frac{1-(\mbox{\boldmath $a$}\cdot\mbox{\boldmath $\nu$}_q)^2}
{\sqrt{\text{det}\,(S_q(\partial B_{d_{\partial D}(p)}(p))-S_q(\partial D))}}.
\end{array}
\tag {2.25}
$$
Here from (2.10) we have
$$\displaystyle
\frac{e^{2\tilde{\tau}d_{\partial D}(p)}}{K(\tau)^2}
=\frac{4\tau^2\epsilon^2}{\eta^2}\frac{e^{2\tilde{\tau}(d_{\partial D}(p)-\eta)}}{(1+O(\tau^{-1}))}.
$$
And also $d_{\partial D}(p)-\eta=\text{dist}\,(D,B)$.
These together with (2.25) yield (2.14).

Next we prove (2.15).
Since we have
$$\displaystyle
\int_D\vert\mbox{\boldmath $V$}'\vert^2dx
=\int_D\vert\nabla\times\mbox{\boldmath $V$}\vert^2dx
+\int_{\partial D}\mbox{\boldmath $\nu$}\cdot\mbox{\boldmath $V$}'\mbox{\boldmath $V$}dS,
\tag {2.26}
$$
from (1.11) we see that it suffices to study the asymptotic behaviour of the second integral on this right-hand
side.
Thus, for this purpose we compute $\mbox{\boldmath $V$}'$.

From (2.7) we have
$$
\displaystyle
(\mbox{\boldmath $M$}(x;p)\mbox{\boldmath $a$})'
=\mbox{\boldmath $a$}\otimes\nabla A-
\left\{\left(\frac{x-p}{\vert x-p\vert}\otimes\frac{x-p}{\vert x-p\vert}\right)\mbox{\boldmath $a$}\right\}\otimes\nabla B
-B\left\{\left(\frac{x-p}{\vert x-p\vert}\otimes\frac{x-p}{\vert x-p\vert}\right)\mbox{\boldmath $a$}\right\}'.
$$
Since
$$\displaystyle
\left(\frac{x-p}{\vert x-p\vert}\right)'
=\frac{1}{\vert x-p\vert}
\left(I_3-\frac{x-p}{\vert x-p\vert}\otimes\frac{x-p}{\vert x-p\vert}\right),
$$
one gets
$$\begin{array}{c}
\displaystyle
\left\{\left(\frac{x-p}{\vert x-p\vert}\otimes\frac{x-p}{\vert x-p\vert}\right)\mbox{\boldmath $a$}\right\}'
=
\left(\frac{x-p}{\vert x-p\vert}\right)'\left(\frac{x-p}{\vert x-p\vert}\cdot\mbox{\boldmath $a$}\right)
+\frac{x-p}{\vert x-p\vert}\otimes\left\{\left(\frac{x-p}{\vert x-p\vert}\right)'\right\}^T\mbox{\boldmath $a$}\\
\\
\displaystyle
=
\frac{1}{\vert x-p\vert}
\left\{
\left(I_3-\frac{x-p}{\vert x-p\vert}\otimes\frac{x-p}{\vert x-p\vert}\right)
\left(\frac{x-p}{\vert x-p\vert}\cdot\mbox{\boldmath $a$}\right)
+\frac{x-p}{\vert x-p\vert}\otimes
\left(I_3-\frac{x-p}{\vert x-p\vert}\otimes\frac{x-p}{\vert x-p\vert}\right)
\mbox{\boldmath $a$}\right\}
\\
\\
\displaystyle
=\frac{1}{\vert x-p\vert}
\left(
\mbox{\boldmath $\omega$}_x\cdot\mbox{\boldmath $a$}
I_3
-2\mbox{\boldmath $\omega$}_x\cdot\mbox{\boldmath $a$}\,\mbox{\boldmath $\omega$}_x
\otimes
\mbox{\boldmath $\omega$}_x
+\mbox{\boldmath $\omega$}_x
\otimes
\mbox{\boldmath $a$}
\right).
\end{array}
$$
Inserting this together with the direct computation results of $\nabla A$ and $\nabla B$ into the expression above
we have
$$
\begin{array}{c}
\displaystyle
(\mbox{\boldmath $M$}(x;p)\mbox{\boldmath $a$})'
=\frac{1}{\tilde{\tau}}\left(\frac{1}{\vert x-p\vert^2}+\frac{2}{\tilde{\tau}\vert x-p\vert^3}\right)
\left(3\mbox{\boldmath $\omega$}_x\otimes\mbox{\boldmath $\omega$}_x\,\mbox{\boldmath $\omega$}_x\cdot\mbox{\boldmath $a$}
-\mbox{\boldmath $a$}\otimes\mbox{\boldmath $\omega$}_x\right)\\
\\
\displaystyle
-
\left\{\frac{1}{\vert x-p\vert}+\frac{3}{\tilde{\tau}}
\left(\frac{1}{\vert x-p\vert^2}+\frac{1}{\tilde{\tau}\vert x-p\vert^3}\right)
\right\}\left(
\mbox{\boldmath $\omega$}_x\cdot\mbox{\boldmath $a$}
I_3
-2\mbox{\boldmath $\omega$}_x\cdot\mbox{\boldmath $a$}\,\mbox{\boldmath $\omega$}_x
\otimes
\mbox{\boldmath $\omega$}_x
+\mbox{\boldmath $\omega$}_x
\otimes
\mbox{\boldmath $a$}
\right).
\end{array}
$$
In particular, we have as $\tau\longrightarrow\infty$,
$$\begin{array}{c}
\displaystyle
(\mbox{\boldmath $M$}(x;p)\mbox{\boldmath $a$})'
=-
\frac{1}{\vert x-p\vert}
\left(
\mbox{\boldmath $\omega$}_x\cdot\mbox{\boldmath $a$}
I_3
-2\mbox{\boldmath $\omega$}_x\cdot\mbox{\boldmath $a$}\,\mbox{\boldmath $\omega$}_x
\otimes
\mbox{\boldmath $\omega$}_x
+\mbox{\boldmath $\omega$}_x
\otimes
\mbox{\boldmath $a$}
\right)
+O\left(\frac{1}{\tau}\right)
\end{array}
\tag {2.27}
$$
uniformly for $x\in\partial D$.

On the other hand, we have
$$\begin{array}{c}
\displaystyle
(\mbox{\boldmath $M$}(x;p)\mbox{\boldmath $a$})\otimes\nabla v(x)
=-\left(\frac{\tilde{\tau}}{\vert x-p\vert}+\frac{1}{\vert x-p\vert^2}\right)v(x)
(\mbox{\boldmath $M$}(x;p)\mbox{\boldmath $a$})\otimes(x-p)\\
\\
\displaystyle
=-\tilde{\tau}v(x)
\left(1+\frac{1}{\tilde{\tau}\vert x-p\vert}\right)
(\mbox{\boldmath $M$}(x;p)\mbox{\boldmath $a$})\otimes\mbox{\boldmath $\omega$}_x
\end{array}
$$
and thus (2.21) gives, as $\tau\longrightarrow\infty$
$$\begin{array}{c}
\displaystyle
(\mbox{\boldmath $M$}(x;p)\mbox{\boldmath $a$})\otimes\nabla v(x)
=-\tilde{\tau}v(x)
\left\{(I_3-\mbox{\boldmath $\omega$}_x\otimes\mbox{\boldmath $\omega$}_x)\mbox{\boldmath $a$}\otimes\mbox{\boldmath $\omega$}_x
+O\left(\frac{1}{\tau}\right)\right\}
\end{array}
\tag {2.28}
$$
uniformly for $x\in\partial D$.

From (2.6) we have
$$\displaystyle
\mbox{\boldmath $V$}'(x)=K(\tau)\tilde{f}(\tau)\left\{v(x)(\mbox{\boldmath $M$}(x;p)\mbox{\boldmath $a$})'
+(\mbox{\boldmath $M$}(x;p)\mbox{\boldmath $a$})\otimes\nabla v(x)\right\}.
$$
This together with (2.27) and (2.28) yields, as $\tau\longrightarrow\infty$
$$\displaystyle
\mbox{\boldmath $V$}'(x)
=-\tilde{\tau}K(\tau)\tilde{f}(\tau)v(x)
\left\{(I_3-\mbox{\boldmath $\omega$}_x\otimes\mbox{\boldmath $\omega$}_x)\mbox{\boldmath $a$}\otimes\mbox{\boldmath $\omega$}_x
+O\left(\frac{1}{\tau}\right)\right\}.
\tag {2.29}
$$
On the other hand, from (2.6) and (2.21) we obtain
$$\displaystyle
\mbox{\boldmath $V$}(x)
=K(\tau)\tilde{f}(\tau)v(x)
\left\{
\left(I_3-\mbox{\boldmath $\omega$}_x
\otimes\mbox{\boldmath $\omega$}_x\right)\mbox{\boldmath $a$}
+O\left(\frac{1}{\tau}\right)
\right\}.
$$
A combination of this and (2.29) gives
$$\begin{array}{c}
\displaystyle
\mbox{\boldmath $V$}'(x)\mbox{\boldmath $V$}(x)\\
\\
\displaystyle
=-\tilde{\tau}K(\tau)^2\tilde{f}(\tau)^2v(x)^2
\left\{(I_3-\mbox{\boldmath $\omega$}_x\otimes\mbox{\boldmath $\omega$}_x)\mbox{\boldmath $a$}\otimes\mbox{\boldmath $\omega$}_x
\left(I_3-
\mbox{\boldmath $\omega$}_x\otimes\mbox{\boldmath $\omega$}_x
\right)\mbox{\boldmath $a$}
+O\left(\frac{1}{\tau}\right)
\right\}.
\end{array}
\tag {2.30}
$$
Since a direct computation yields
$$\begin{array}{c}
\displaystyle
(I_3-\mbox{\boldmath $\omega$}_x\otimes\mbox{\boldmath $\omega$}_x)\mbox{\boldmath $a$}\otimes\mbox{\boldmath $\omega$}_x
\left(I_3-
\mbox{\boldmath $\omega$}_x\otimes\mbox{\boldmath $\omega$}_x
\right)\mbox{\boldmath $a$}
=\mbox{\boldmath $0$},
\end{array}
$$
it follows from (2.30) that
$$\begin{array}{c}
\displaystyle
\vert\mbox{\boldmath $V$}'(x)\mbox{\boldmath $V$}(x)\vert
\le CK(\tau)^2\tilde{f}(\tau)^2v(x)^2.
\end{array}
$$
Thus we obtain
$$\displaystyle
\frac{1}{K(\tau)^2\tilde{f}(\tau)^2}
\left\vert\int_{\partial D}\mbox{\boldmath $\nu$}\cdot\mbox{\boldmath $V$}'\mbox{\boldmath $V$}dS\right\vert
\le
C\int_{\partial D}v^2dS.
$$
Write
$$\begin{array}{c}
\displaystyle
\frac{\displaystyle\left\vert\int_{\partial D}\mbox{\boldmath $\nu$}\cdot\mbox{\boldmath $V$}'\mbox{\boldmath $V$}dS\right\vert}
{\displaystyle J(\tau)}
\le
\frac{1}{\tilde{\tau}}
\frac{\displaystyle C\tilde{\tau}e^{2\tilde{\tau}d_{\partial D}(p)}\int_{\partial D}v^2dS}
{\displaystyle\frac{e^{2\tilde{\tau}d_{\partial D}(p)}}{K(\tau)^2\tilde{f}(\tau)^2}J(\tau)}.
\end{array}
$$
Then applying (2.24) and (2.25) to this right-hand side, we conclude
$$\displaystyle
\int_{\partial D}\mbox{\boldmath $\nu$}\cdot\mbox{\boldmath $V$}'\mbox{\boldmath $V$}dS
=O\left(\frac{J(\tau)}{\tau}\right).
$$
Now a combination of this and (2.26) yields (2.15).

\noindent
$\Box$

\section{Proof of Theorem 1.1}

The proof of Theorem 1.1 starts with establishing the following
brief asymptotic formula of the indicator function.

\proclaim{\noindent Proposition 3.1.}
It holds that, as $\tau\longrightarrow\infty$
$$\begin{array}{c}
\displaystyle
\int_{\Bbb R^3\setminus\overline D}\mbox{\boldmath $f$}(x,\tau)\cdot(\mbox{\boldmath $W$}_e-\mbox{\boldmath $V$})dx
=J(\tau)+E(\tau)+O(\tau^{-3/2}e^{-\tau T}).
\end{array}
\tag {3.1}
$$

\endproclaim

{\it\noindent Proof.} Set $\mbox{\boldmath $R$}=\mbox{\boldmath $W$}_e-\mbox{\boldmath $V$}$.
The proof is divided into two steps.

\noindent
{\it Step 1.}
First we show that
$$\begin{array}{c}
\displaystyle
\int_{\Bbb R^3\setminus\overline D}\mbox{\boldmath $f$}(x,\tau)\cdot\mbox{\boldmath $R$}\,dx
=J(\tau)+E(\tau)\\
\\
\displaystyle
-e^{-\tau T}\left(\int_{\Bbb R^3\setminus\overline D}\mbox{\boldmath $F$}(x,\tau)\cdot\mbox{\boldmath $R$}\,dx
-\int_{\Bbb R^3\setminus\overline D}\mbox{\boldmath $F$}(x,\tau)\cdot\mbox{\boldmath $V$}dx\right),
\end{array}
\tag {3.2}
$$
where
$$\displaystyle
\mbox{\boldmath $F$}(x,\tau)=-\left(\tau\mbox{\boldmath $E$}(x,T)+\frac{1}{\epsilon}\nabla\times\mbox{\boldmath $H$}(x,T)\right)
\tag {3.3}
$$
and $(\mbox{\boldmath $E$},\mbox{\boldmath $H$})$ is the solution of (1.1).

Define
$$\begin{array}{c}
\displaystyle
\mbox{\boldmath $W$}_m(x,\tau)
=\int_0^Te^{-\tau t}\mbox{\boldmath $H$}(x,t)dt,\,\,x\in\Bbb R^3\setminus\overline D.
\end{array}
$$
It is easy to see that integration by parts yields
$$
\displaystyle
\nabla\times\mbox{\boldmath $W$}_e+\tau\mu \mbox{\boldmath $W$}_m=-e^{-\tau T}\mu\mbox{\boldmath $H$}(x,T)
\,\,\text{in}\,\Bbb R^3\setminus\overline D,
\tag {3.4}
$$
$$
\displaystyle
\nabla\times\mbox{\boldmath $W$}_m-\tau\epsilon \mbox{\boldmath $W$}_e
-\frac{\epsilon}{\tau}\mbox{\boldmath $f$}(x,\tau)
=e^{-\tau T}\epsilon \mbox{\boldmath $E$}(x,T)\,\,\text{in}\,\Bbb R^3\setminus\overline D
\tag {3.5}
$$
and
$$
\displaystyle
\mbox{\boldmath $\nu$}\times\mbox{\boldmath $W$}_e=\mbox{\boldmath $0$}\,\,\text{on}\,\partial D.
\tag {3.6}
$$
Taking the rotation of (3.4) and (3.5), respectively, we obtain the following equation:
$$\displaystyle
\frac{1}{\mu\epsilon}\nabla\times\nabla\times\mbox{\boldmath $W$}_e
+\tau^2\mbox{\boldmath $W$}_e+\mbox{\boldmath $f$}(x,\tau)
=e^{-\tau T}\mbox{\boldmath $F$}(x,\tau)
\,\,\text{in}\,\Bbb R^3\setminus\overline D.
\tag {3.7}
$$
Integration by parts gives
$$\begin{array}{c}
\displaystyle
\int_{\Bbb R^3\setminus\overline D}\left\{(\nabla\times\nabla\times\mbox{\boldmath $W$}_e)\cdot\mbox{\boldmath $V$}
-(\nabla\times\nabla\times\mbox{\boldmath $V$})\cdot\mbox{\boldmath $W$}_e\right\}dx\\
\\
\displaystyle
=\int_{\partial D}
\left\{(\mbox{\boldmath $\nu$}\times(\nabla\times\mbox{\boldmath $V$}))\cdot\mbox{\boldmath $W$}_e
-(\mbox{\boldmath $\nu$}\times(\nabla\times\mbox{\boldmath $W$}_e))\cdot\mbox{\boldmath $V$}
\right\}dS.
\end{array}
$$
(3.6) ensures that the first term on this right-hand side vanishes.
And we have
$$
\displaystyle
(\mbox{\boldmath $\nu$}\times(\nabla\times\mbox{\boldmath $W$}_e))\cdot\mbox{\boldmath $V$}
=(\nabla\times\mbox{\boldmath $W$}_e)\times\mbox{\boldmath $V$}\cdot\mbox{\boldmath $\nu$}
=(\mbox{\boldmath $V$}\times\mbox{\boldmath $\nu$})\cdot
(\nabla\times\mbox{\boldmath $W$}_e)
=-(\mbox{\boldmath $\nu$}\times\mbox{\boldmath $V$})\cdot
(\nabla\times\mbox{\boldmath $W$}_e).
$$
Thus
$$\begin{array}{c}
\displaystyle
\int_{\Bbb R^3\setminus\overline D}\left\{(\nabla\times\nabla\times\mbox{\boldmath $W$}_e)\cdot\mbox{\boldmath $V$}
-(\nabla\times\nabla\times\mbox{\boldmath $V$})\cdot\mbox{\boldmath $W$}_e\right\}dx
=\int_{\partial D}(\mbox{\boldmath $\nu$}\times\mbox{\boldmath $V$})\cdot
(\nabla\times\mbox{\boldmath $W$}_e)dS.
\end{array}
$$
Substituting (1.3) and (3.7) into this, we obtain
$$\begin{array}{c}
\displaystyle
\frac{1}{\mu\epsilon}\int_{\partial D}(\mbox{\boldmath $\nu$}\times\mbox{\boldmath $V$})\cdot\nabla\times\mbox{\boldmath $W$}_e\,dS\\
\\
\displaystyle
=\int_{\Bbb R^3\setminus\overline D}\mbox{\boldmath $f$}(x,\tau)\cdot\mbox{\boldmath $R$}\,dx
+e^{-\tau T}\int_{\Bbb R^3\setminus\overline D}\mbox{\boldmath $F$}(x,\tau)\cdot\mbox{\boldmath $V$}dx.
\end{array}
\tag {3.8}
$$
Write
$$\begin{array}{c}
\displaystyle
\frac{1}{\mu\epsilon}\int_{\partial D}(\mbox{\boldmath $\nu$}\times\mbox{\boldmath $V$})\cdot\nabla\times\mbox{\boldmath $W$}_e\,dS\\
\\
\displaystyle
=\frac{1}{\mu\epsilon}\int_{\partial D}(\mbox{\boldmath $\nu$}\times\mbox{\boldmath $V$})\cdot\nabla\times\mbox{\boldmath $V$}\,dS
+\frac{1}{\mu\epsilon}\int_{\partial D}(\mbox{\boldmath $\nu$}\times\mbox{\boldmath $V$})\cdot\nabla\times\mbox{\boldmath $R$}\,dS.
\end{array}
\tag {3.9}
$$
Since $\mbox{\boldmath $R$}$ satisfies
$$
\displaystyle
\frac{1}{\mu\epsilon}\nabla\times\nabla\times\mbox{\boldmath $R$}+\tau^2
\mbox{\boldmath $R$}=e^{-\tau T}\mbox{\boldmath $F$}(x,\tau)
\,\,\text{in}\,\Bbb R^3\setminus\overline D
\tag {3.10}
$$
and
$$
\displaystyle
\mbox{\boldmath $\nu$}\times\mbox{\boldmath $R$}=-\mbox{\boldmath $\nu$}\times\mbox{\boldmath $V$}
\,\,\text{on}\,\partial D,
\tag {3.11}
$$
integration by parts gives
$$\begin{array}{c}
\displaystyle
e^{-\tau T}\int_{\Bbb R^3\setminus\overline D}\mbox{\boldmath $F$}(x,\tau)\cdot
\mbox{\boldmath $R$}\,dx
=-\frac{1}{\mu\epsilon}\int_{\partial D}(\mbox{\boldmath $\nu$}\times\mbox{\boldmath $V$})
\cdot\nabla\times
\mbox{\boldmath $R$}\,dS+E(\tau),
\end{array}
$$
that is,
$$\begin{array}{c}
\displaystyle
\frac{1}{\mu\epsilon}\int_{\partial D}(\mbox{\boldmath $\nu$}\times\mbox{\boldmath $V$})
\cdot\nabla\times
\mbox{\boldmath $R$}\,dS
=E(\tau)
-e^{-\tau T}\int_{\Bbb R^3\setminus\overline D}\mbox{\boldmath $F$}(x,\tau)\cdot
\mbox{\boldmath $R$}\,dx.
\end{array}
\tag {3.12}
$$
Now (3.2) follows from (1.11), (1.12), (3.8), (3.9), (3.12) and (2.16).

\noindent
{\it Step 2.}
It follows from the definition of the weak solution of (1.3) that
$$\displaystyle
\frac{1}{\mu\epsilon}\int_{\Bbb R^3}\vert\nabla\times\mbox{\boldmath $V$}\vert^2dx
+\tau^2\int_{\Bbb R^3}\left\vert\mbox{\boldmath $V$}+\frac{\mbox{\boldmath $f$}}{2\tau^2}\right\vert^2dx
=\frac{1}{4\tau^2}\int_{\Bbb R^3}\vert\mbox{\boldmath $f$}\vert^2dx.
\tag {3.13}
$$
Since $\tilde{f}(\tau)=O(\tau^{-3/2})$, we have
$$\displaystyle
\Vert\mbox{\boldmath $f$}(\,\cdot\,,\tau)\Vert_{L^2(\Bbb R^3)}
=O(\tau^{-1/2}).
\tag {3.14}
$$
Then, applying  the inequality
$$\displaystyle
\vert\mbox{\boldmath $A$}+\mbox{\boldmath $B$}\vert^2
\ge\frac{1}{2}\vert\mbox{\boldmath $A$}\vert^2-\vert\mbox{\boldmath $B$}\vert^2
\tag {3.15}
$$
to the second term in the left-hand side on (3.13), we obtain, as
$\tau\longrightarrow\infty$
$$\displaystyle
\frac{1}{\mu\epsilon}\int_{\Bbb R^3}\vert\nabla\times\mbox{\boldmath $V$}\vert^2dx
+\tau^2\int_{\Bbb R^3}\left\vert\mbox{\boldmath $V$}\right\vert^2dx
=O(\tau^{-3}).
\tag {3.16}
$$
%$$\displaystyle
%J(\tau)=O(\tau^{-3}).
%\tag {3.16}
%$$
Next we prove that, as $\tau\longrightarrow\infty$
$$\displaystyle
E(\tau)=O(\tau^{-3}).
\tag {3.17}
$$
Write
$$\begin{array}{c}
\displaystyle
\tau^2\vert\mbox{\boldmath $R$}\vert^2-\mbox{\boldmath $f$}\cdot\mbox{\boldmath $R$}
-e^{-\tau T}\mbox{\boldmath $F$}\cdot\mbox{\boldmath $R$}
=\tau^2\left\vert\mbox{\boldmath $R$}
-\frac{\mbox{\boldmath $f$}+e^{-\tau T}\mbox{\boldmath $F$}}{2\tau^2}\right\vert^2
-\frac{\vert\mbox{\boldmath $f$}+e^{-\tau T}\mbox{\boldmath $F$}\vert^2}{4\tau^2}.
\end{array}
$$
Substituting this into (3.2), we obtain
$$\begin{array}{c}
\displaystyle
\frac{1}{\mu\epsilon}\int_{\Bbb R^3\setminus\overline D}\vert\nabla\times\mbox{\boldmath $R$}\vert^2dx
+\tau^2\int_{\Bbb R^3\setminus\overline D}\left\vert\mbox{\boldmath $R$}
-\frac{\mbox{\boldmath $f$}+e^{-\tau T}\mbox{\boldmath $F$}}{2\tau^2}\right\vert^2dx
+J(\tau)
\\
\\
\displaystyle
=\frac{1}{4\tau^2}\int_{\Bbb R^3\setminus\overline D}\vert\mbox{\boldmath $f$}+e^{-\tau T}\mbox{\boldmath $F$}\vert^2dx
+e^{-\tau T}\int_{\Bbb R^3\setminus\overline D}\mbox{\boldmath $F$}\cdot\mbox{\boldmath $V$}dx.
\end{array}
\tag {3.18}
$$
Dropping the second and third terms on the left-hand side of (3.18), we obtain
$$\begin{array}{c}
\displaystyle
\frac{1}{\mu\epsilon}\int_{\Bbb R^3\setminus\overline D}\vert\nabla\times\mbox{\boldmath $R$}\vert^2dx
\le
\frac{1}{4\tau^2}
\int_{\Bbb R^3\setminus\overline D}\vert\mbox{\boldmath $f$}+e^{-\tau T}\mbox{\boldmath $F$}\vert^2dx
+e^{-\tau T}\int_{\Bbb R^3\setminus\overline D}\vert\mbox{\boldmath $F$}\cdot\mbox{\boldmath $V$}\vert dx.
\end{array}
\tag {3.19}
$$

By (3.3) we have
$$\displaystyle
\Vert\mbox{\boldmath $F$}\Vert_{L^2(\Bbb R^3\setminus\overline D)}=O(\tau).
\tag {3.20}
$$
It follows from (3.16) that
$$\displaystyle
\Vert \mbox{\boldmath $V$}\Vert_{L^2(\Bbb R^3\setminus\overline D)}=O(\tau^{-5/2}).
$$
Applying these and (3.14) to the right-hand side on (3.19), we obtain
$$\displaystyle
\frac{1}{\mu\epsilon}\int_{\Bbb R^3\setminus\overline D}\vert\nabla\times\mbox{\boldmath $R$}\vert^2dx
=O(\tau^{-3}).
$$

On the other hand, dropping the first and third terms on the left-hand side on (3.18) and using (3.15), we obtain
$$\begin{array}{c}
\displaystyle
\frac{\tau^2}{2}\int_{\Bbb R^3\setminus\overline D}\left\vert\mbox{\boldmath $R$}\right\vert^2dx
\le\frac{1}{2\tau^2}\int_{\Bbb R^3\setminus\overline D}\vert\mbox{\boldmath $f$}+e^{-\tau T}\mbox{\boldmath $F$}\vert^2dx
+e^{-\tau T}\int_{\Bbb R^3\setminus\overline D}\vert\mbox{\boldmath $F$}\cdot\mbox{\boldmath $V$}\vert dx.
\end{array}
$$
Thus, by the same reason above we obtain
$$\displaystyle
\tau^2\int_{\Bbb R^3\setminus\overline D}\left\vert\mbox{\boldmath $R$}\right\vert^2dx
=O(\tau^{-3}).
$$
This completes the proof of (3.17).

Finally from (3.16), (3.17) and (3.20) we have
$$\displaystyle
\int_{\Bbb R^3\setminus\overline D}\mbox{\boldmath $F$}(x,\tau)\cdot\mbox{\boldmath $R$}\,dx
+\int_{\Bbb R^3\setminus\overline D}\mbox{\boldmath $F$}(x,\tau)\cdot\mbox{\boldmath $V$}dx
=O(\tau^{-3/2}).
$$
Thus, a combination of this and (3.2) yields (3.1).

\noindent
$\Box$

{\bf\noindent Remark 3.1.}
From (3.4) and (3.5) we obtain also the following equation for $\mbox{\boldmath $W$}_m$:
$$\displaystyle
\frac{1}{\mu\epsilon}\nabla\times\nabla\times\mbox{\boldmath $W$}_m
+\tau^2\mbox{\boldmath $W$}_m
-\frac{\epsilon}{\tau}\nabla\times\mbox{\boldmath $f$}(x,\tau)
=e^{-\tau T}
\tilde{\mbox{\boldmath $F$}}(x,\tau)
\,\,\text{in}\,\Bbb R^3\setminus\overline D,
$$
where
$$\displaystyle
\tilde{\mbox{\boldmath $F$}}(x,\tau)=\frac{1}{\mu}\nabla\times\mbox{\boldmath $E$}(x,T)-\tau\mbox{\boldmath $H$}(x,T).
$$
In this paper, we will not make use of this equation.

Next we prove

\proclaim{\noindent Lemma 3.1.}
There exist positive constant $C$ and $\tau_0$ such that, for all $\tau\ge\tau_0$
$$\displaystyle
E(\tau)\le C(\tau^2J(\tau)+e^{-2\tau T}).
\tag {3.21}
$$

\endproclaim

{\it\noindent Proof.}
Set $\mbox{\boldmath $R$}=\mbox{\boldmath $W$}_e-\mbox{\boldmath $V$}$.
Taking the scalar product of equation (3.10) with $\mbox{\boldmath $R$}$,
integrating over $\Bbb R^3\setminus\overline D$ and using boundary condition (3.11) on $\partial D$ and (1.12), we have
$$\begin{array}{c}
\displaystyle
E(\tau)
=\frac{1}{\mu\epsilon}\int_{\partial D}\mbox{\boldmath $\nu$}\times\mbox{\boldmath $V$}
\cdot\nabla\times\mbox{\boldmath $R$}dS
+e^{-\tau T}\int_{\Bbb R^3\setminus\overline D}\mbox{\boldmath $F$}\cdot\mbox{\boldmath $R$}dx.
\end{array}
\tag {3.22}
$$
By the trace theorem (\cite{N}, p. 209, Theorem 5.4.2.), one can choose a lifting $\tilde{\mbox{\boldmath $V$}}$ of
$\mbox{\boldmath $\nu$}\times\mbox{\boldmath $V$}$ on $\partial D$ in such a way that
$$
\displaystyle
\Vert\tilde{\mbox{\boldmath $V$}}\Vert_{L^2(\Bbb R^3\setminus\overline D)}^2
+\Vert\nabla\times\tilde{\mbox{\boldmath $V$}}\Vert_{L^2(\Bbb R^3\setminus\overline D)}^2
\le C^2\Vert\mbox{\boldmath $\nu$}\times\mbox{\boldmath $V$}\Vert_{H^{-1/2}_{\text{div}}(\partial D)}^2.
$$
See also p. 191 in \cite{N} for the definition of $H^{-1/2}_{\text{div}}(\partial D)$ and the norm.
Note that $C$ is a positive constant and independent of $\mbox{\boldmath $V$}$.

Again the trace theorem tells us also that
$$\displaystyle
\Vert\mbox{\boldmath $\nu$}\times\mbox{\boldmath $V$}\Vert_{H^{-1/2}_{\text{div}}(\partial D)}^2
\le(C')^2(\Vert\mbox{\boldmath $V$}\Vert_{L^2(D)}^2
+\Vert\nabla\times\mbox{\boldmath $V$}\Vert_{L^2(D)}^2),
$$
where $C'$ is a positive constant and independent of $\mbox{\boldmath $V$}$.
Thus, we have
$$\displaystyle
\Vert\tilde{\mbox{\boldmath $V$}}\Vert_{L^2(\Bbb R^3\setminus\overline D)}^2
+\Vert\nabla\times\tilde{\mbox{\boldmath $V$}}\Vert_{L^2(\Bbb R^3\setminus\overline D)}^2
\le (CC')^2(\Vert\mbox{\boldmath $V$}\Vert_{L^2(D)}^2
+\Vert\nabla\times\mbox{\boldmath $V$}\Vert_{L^2(D)}^2).
\tag {3.23}
$$
Moreover, from equation (3.10) one gets
$$\begin{array}{c}
\displaystyle
\frac{1}{\mu\epsilon}\int_{\partial D}\mbox{\boldmath $\nu$}\times\mbox{\boldmath $V$}\cdot\nabla\times\mbox{\boldmath $R$}dS\\
\\
\displaystyle
=-\frac{1}{\mu\epsilon}
\int_{\Bbb R^3\setminus\overline D}\nabla\times\mbox{\boldmath $R$}\cdot\nabla\times\tilde{\mbox{\boldmath $V$}}dx
-\tau^2\int_{\Bbb R^3\setminus\overline D}\mbox{\boldmath $R$}\cdot\tilde{\mbox{\boldmath $V$}}dx
+e^{-\tau T}\int_{\Bbb R^3\setminus\overline D}\mbox{\boldmath $F$}\cdot\tilde{\mbox{\boldmath $V$}}dx.
\end{array}
$$
Substituting this into (3.22) and estimating from above, we obtain
$$\begin{array}{c}
\displaystyle
E(\tau)
\le
\frac{1}{\mu\epsilon}
\Vert\nabla\times\mbox{\boldmath $R$}\Vert_{L^2(\Bbb R^3\setminus\overline D)}
\Vert\nabla\times\tilde{\mbox{\boldmath $V$}}\Vert_{L^2(\Bbb R^3\setminus\overline D)}
+\tau^2\Vert\mbox{\boldmath $R$}\Vert_{L^2(\Bbb R^3\setminus\overline D)}
\Vert\tilde{\mbox{\boldmath $V$}}\Vert_{L^2(\Bbb R^3\setminus\overline D)}\\
\\
\displaystyle
+e^{-\tau T}\Vert\mbox{\boldmath $F$}\Vert_{L^2(\Bbb R^3\setminus\overline D)}
(\Vert\tilde{\mbox{\boldmath $V$}}\Vert_{L^2(\Bbb R^3\setminus\overline D)}
+
\Vert\mbox{\boldmath $R$}\Vert_{L^2(\Bbb R^3\setminus\overline D)}).
\end{array}
\tag {3.24}
$$
Here we make use of the following trivial estimates
$$
\displaystyle
\Vert\nabla\times\mbox{\boldmath $R$}\Vert_{L^2(\Bbb R^3\setminus\overline D)}
\le\sqrt{\mu\epsilon}\sqrt{E(\tau)}
\tag {3.25}
$$
and
$$
\displaystyle
\Vert\mbox{\boldmath $R$}\Vert_{L^2(\Bbb R^3\setminus\overline D)}
\le\tau^{-1}\sqrt{E(\tau)}.
\tag {3.26}
$$
Applying (3.20), (3.25) and (3.26) to the right-hand side on (3.24), we obtain
$$\begin{array}{c}
\displaystyle
E(\tau)
\le
C_1\left(\Vert\nabla\times\tilde{\mbox{\boldmath $V$}}\Vert_{L^2(\Bbb R^3\setminus\overline D)}
+\tau\Vert\tilde{\mbox{\boldmath $V$}}\Vert_{L^2(\Bbb R^3\setminus\overline D)}
\right)
\sqrt{E(\tau)}\\
\\
\displaystyle
+C_2e^{-\tau T}\tau\Vert\tilde{\mbox{\boldmath $V$}}\Vert_{L^2(\Bbb R^3\setminus\overline D)}
+C_3e^{-\tau T}\sqrt{E(\tau)}.
\end{array}
$$
Then, applying a standard technique to the first and last terms on this right-hand side,
we obtain
$$
\begin{array}{c}
\displaystyle
E(\tau)
\le
C_4\left(\Vert\nabla\times\tilde{\mbox{\boldmath $V$}}\Vert_{L^2(\Bbb R^3\setminus\overline D)}^2
+\tau^2\Vert\tilde{\mbox{\boldmath $V$}}\Vert_{L^2(\Bbb R^3\setminus\overline D)}^2
\right)
+C_5e^{-2\tau T}.
\end{array}
$$
Then, a combination of this, (3.23) and trivial inequality
$$\displaystyle
\Vert\mbox{\boldmath $V$}\Vert_{L^2(D)}^2
+\Vert\nabla\times\mbox{\boldmath $V$}\Vert_{L^2(D)}^2
\le C(1+\tau^{-2})J(\tau)
$$
yields (3.21).

\noindent
$\Box$

Now from (3.1), (2.8), (2.9) and (3.21) we obtain
$$\displaystyle
\limsup_{\tau\longrightarrow\infty}\tau e^{\displaystyle 2\tau\sqrt{\mu\epsilon}\text{dist}\,(D,B)}
\int_{\Bbb R^3\setminus\overline D}
\mbox{\boldmath $f$}(x,\tau)\cdot(\mbox{\boldmath $W$}_e-\mbox{\boldmath $V$})dx<\infty
$$
and
$$\displaystyle
\liminf_{\tau\longrightarrow\infty}\tau^{5+2\gamma}e^{\displaystyle 2\tau\sqrt{\mu\epsilon}\text{dist}\,(D,B)}
\int_{\Bbb R^3\setminus\overline D}
\mbox{\boldmath $f$}(x,\tau)\cdot(\mbox{\boldmath $W$}_e-\mbox{\boldmath $V$})dx>0
$$
provided $T>2\sqrt{\mu\epsilon}\text{dist}\,(D,B)$.  From these we immediately obtain Theorem 1.1.

{\bf\noindent Remark 3.2.}
(3.21) in Lemma 3.1 is {\it not sharp}, however, for Theorem 1.1 it is enough.
For Theorem 1.2 we need more accurate estimate like $E(\tau)\sim J(\tau)$.

\section{Proof of Theorem 1.2}

First it is easy to see that Theorem 1.2 follows from the following theorem and (2.14).

\proclaim{\noindent Theorem 4.1.}
Assume that $\Lambda_{\partial D}(p)$ is finite and that (1.8) and (1.9) are satisfied.
Let $f$ satisfy (1.6) for a $\gamma\in\Bbb R$.
Let $T>2\sqrt{\mu\epsilon}\text{dist}\,(D,B)$.  Then,
as $\tau\longrightarrow\infty$, we have
$$\begin{array}{c}
\displaystyle
\int_{\Bbb R^3\setminus\overline D}\mbox{\boldmath $f$}(x,\tau)\cdot(\mbox{\boldmath $W$}_e-\mbox{\boldmath $V$})dx
=2J(\tau)(1+O(\tau^{-1/2})).
\end{array}
\tag {4.1}
$$
\endproclaim

Thus the purpose of this section is to describe the proof of
Theorem 4.1. However, note that under the assumption (1.6) for a
$\gamma\in\Bbb R$ it holds that
$$\displaystyle
\frac{\displaystyle
\tau^2 e^{2\tau\sqrt{\mu\epsilon}\text{dist}\,(D,B)}}
{\displaystyle\tilde{f}(\tau)^2}
\tau^{-1/2}e^{-\tau T}
=O(\tau^{2-1/2+2\gamma}e^{-\tau(T-2\sqrt{\mu\epsilon}\text{dist}\,(D,B)}).
$$
Thus, if we have the estimate
$$\displaystyle
E(\tau)=J(\tau)(1+O(\tau^{-1/2})),
\tag {4.2}
$$
then (3.1) yields (4.1).  Thus, the proof of (4.1) is reduced to that of (4.2)
which shows the asymptotic coincidence of $E(\tau)$ and $J(\tau)$ as $\tau\longrightarrow\infty$.

The proof of (4.2) employs the Lax-Phillips reflection
argument in \cite{LP}, however, some technical parts are
different. Anyway that is based on: a representation formula of
$E(\tau)-J(\tau)$ via a reflection. Thus, the following subsection
starts with describing a reflection principle across $\partial D$
from inside to outside.

\subsection{Reflection principle}

One can choose a positive number $\delta_0$ in such a way that: given $x\in\Bbb R^3\setminus D$($x\in\overline D$)
with $d_{\partial D}(x)<2\delta_0$ there exists a unique $q=q(x)\in\partial D$ such that
$x=q+d_{\partial D}(x)\mbox{\boldmath $\nu$}_q$($x=q-d_{\partial D}(x)\mbox{\boldmath $\nu$}_q$).
Both $d_{\partial D}(x)$ and $q(x)$ are $C^k$ therein provided $\partial D$
is $C^k$ with $k\ge 2$.  See Lemma 14.16 in \cite{GT} for this.

For $x$ with $d_{\partial D}(x)<2\delta_0$ define $x^r=2q(x)-x$, $\pi(x)=\mbox{\boldmath $\nu$}_{q(x)}
\otimes\mbox{\boldmath $\nu$}_{q(x)}$ and $\displaystyle
\mbox{\boldmath $n$}(x)=\mbox{\boldmath $\nu$}_{q(x)}$.
Note that $\mbox{\boldmath $n$}$ is $C^3$ if $\partial D$ is $C^4$.

The reflection principle what we say in this paper consists of two parts
summarized as the following propositions.

\proclaim{\noindent Proposition 4.1.}
Assume that $\partial D$ is $C^4$.
Let $\mbox{\boldmath $V$}$ be a vector field over $D$ and $C^2$ in $D$.
For $x\in\Bbb R^3\setminus D$ with $d_{\partial D}(x)<2\delta_0$,
define
$$\displaystyle
\mbox{\boldmath $V$}^*(x)=-\mbox{\boldmath $A$}(x^r)+\mbox{\boldmath $B$}(x^r)+2d_{\partial D}(x)\mbox{\boldmath $n$}'(x)\mbox{\boldmath $A$}(x^r),
\tag {4.3}
$$
where $\mbox{\boldmath $A$}(y)=(I-\pi(y))\mbox{\boldmath $V$}(y)$ and $\mbox{\boldmath $B$}(y)=\pi(y)\mbox{\boldmath $V$}(y)$
for $y\in D$ with $d_{\partial D}(y)<2\delta_0$.

Then, $\mbox{\boldmath $V$}^*$ satisfies
$$
\displaystyle
\mbox{\boldmath $V$}^*\times\mbox{\boldmath $\nu$}=-\mbox{\boldmath $V$}\times\mbox{\boldmath $\nu$}\,\,\,\text{on}\,\partial D
\tag {4.4}
$$
and
$$
\displaystyle
\mbox{\boldmath $\nu$}\times(\nabla\times\mbox{\boldmath $V$}^*)=\mbox{\boldmath $\nu$}\times(\nabla\times\mbox{\boldmath $V$})\,\,\,\text{on}\,\partial D.
\tag {4.5}
$$

\endproclaim

{\it\noindent Proof.}
Define
$$\displaystyle
\tilde{\mbox{\boldmath $V$}}(x)=-\mbox{\boldmath $A$}(x^r)+\mbox{\boldmath $B$}(x^r)
\tag {4.6}
$$
and
$$\displaystyle
\mbox{\boldmath $C$}(x)=2d_{\partial D}(x)\mbox{\boldmath $n$}'(x)\mbox{\boldmath $V$}(x^r).
\tag {4.7}
$$
We have
$$\displaystyle
\mbox{\boldmath $V$}^*(x)=\tilde{\mbox{\boldmath $V$}}(x)+\mbox{\boldmath $C$}(x).
$$

First, we claim that $\tilde{\mbox{\boldmath $V$}}$ satisfies the following boundary conditions.

{\bf\noindent Claim 1.}
$\tilde{\mbox{\boldmath $V$}}$ satisfies
the following boundary conditions
$$
\displaystyle
\tilde{\mbox{\boldmath $V$}}\times\mbox{\boldmath $\nu$}_x=-\mbox{\boldmath $V$}\times\mbox{\boldmath $\nu$}_x\,\,\text{on}\,\partial D;
\tag {4.8}
$$
$$
\displaystyle
\mbox{\boldmath $\nu$}_x\times(\nabla\times\tilde{\mbox{\boldmath $V$}})
=\mbox{\boldmath $\nu$}_x\times(\nabla\times{\mbox{\boldmath $V$}})
-2S_x(\partial D)\mbox{\boldmath $A$}\,\,\text{on}\,\partial D,
\tag {4.9}
$$
where $S_x(\partial D)$ denotes the {\it shape operator} of $\partial D$ at $x\in\partial D$ with respect to $\mbox{\boldmath $\nu$}_x$.

Next we claim

{\bf\noindent Claim 2.}
We have
$$\displaystyle
\mbox{\boldmath $\nu$}_x\times(\nabla\times\mbox{\boldmath $C$})=2S_x(\partial D)\mbox{\boldmath $A$}\,\,\text{on}\,\partial D.
\tag {4.10}
$$

Now from trivial identity $\mbox{\boldmath $C$}=\mbox{\boldmath $0$}$ on $\partial D$ and (4.8) we obtain
(4.4); from (4.9) and (4.10) we obtain (4.5).  This completes the proof of Proposition 4.1.
See also Appendix for the proof of Claims 1 and 2.

\noindent
$\Box$

Since the proof of the following proposition is tedious and so is
described in Appendix.

\proclaim{\noindent Proposition 4.2.}
Assume that $\partial D$ is $C^4$.
If $\mbox{\boldmath $V$}$ satisfies
$$\displaystyle
\frac{1}{\mu\epsilon}\nabla\times\nabla\times\mbox{\boldmath $V$}+\tau^2\mbox{\boldmath $V$}=\mbox{\boldmath $0$}\,\,\text{in}\,D,
\tag {4.11}
$$
then, $\mbox{\boldmath $V$}^*$ defined as (4.3) satisfies
$$\begin{array}{c}
\displaystyle
\frac{1}{\mu\epsilon}\nabla\times\nabla\times\mbox{\boldmath $V$}^*+\tau^2\mbox{\boldmath $V$}^*\\
\\
\displaystyle
=\text{terms from $\mbox{\boldmath $V$}(x^r)$ and $\mbox{\boldmath $V$}'(x^r)$}
+2d_{\partial D}(x)
\times\text{terms from $(\nabla^2\mbox{\boldmath $V$})(x^r)$}
\end{array}
\tag {4.12}
$$
and all the coefficients in this right-hand side are {\it
independent} of $\tau$ and continuous, in particular, the
coefficients come from the second order terms are $C^1$ in a
tubular neighbourhood of $\partial D$.

\endproclaim

{\bf\noindent Remark 4.1.}
Note that, for $y\in D$ with $d_{\partial D}(y)<2\delta_0$ we have the decomposition
$$\displaystyle
\mbox{\boldmath $V$}(y)=\mbox{\boldmath $A$}(y)+\mbox{\boldmath $B$}(y).
\tag {4.13}
$$
If $\partial D$ is a {\it plane}, then $\mbox{\boldmath $n$}'(x)\equiv 0$
and the third term in the right-hand side on (4.3) vanishes.
Thus, in this case Propositions 4.1 and 4.2 become
the reflection principle used in \cite{LYZ} for inverse obstacle scattering
for Maxwell's equations in a frequency domain (replaced $\tau^2$ with $-k^2$) .
They employed this principle for a different purpose from us, more precisely,
establishing a uniqueness theorem for {\it polygonal obstacles} in a {\it single frequency domain}.
In the curved boundary case, $\mbox{\boldmath $n'$}\not=\mbox{\boldmath $0$}$ and we need the {\it correction term}
$2\displaystyle d_{\partial D}(x)\mbox{\boldmath $n$}'(x)\mbox{\boldmath $A$}(x^r)$.
For more detailed information about (4.12) see Lemmas A.1 and A.2.

\subsection{Proof of the estimate (4.2)}

In this subsection we start with describing a representation
formula of $E(\tau)-J(\tau)$ in terms of the reflection across
$\partial D$.

Let $0<\delta<\delta_0/2$.  Choose a smooth function $\phi=\phi_{\delta}$ defined on the whole space in such a way that
(i)  $0\le\phi\le 1$; (ii) $\phi(x)=1$ if $d_{\partial D}(x)<\delta$ and $\phi(x)=0$ if $d_{\partial D}(x)>2\delta$;
(iii) $\vert\nabla\phi(x)\vert\le C\delta^{-1}$;
$\vert(\partial^2/\partial x_i\partial x_j)\phi(x)\vert\le C\delta^{-2}$ with $i,j=1,\cdots,3$.

Define
$$\displaystyle
\mbox{\boldmath $V$}^r(x)=\phi(x)\mbox{\boldmath $V$}^*(x),\,\,x\in\Bbb R^3,
\tag {4.14}
$$
where $\mbox{\boldmath $V$}^*$ is given by (4.3).

Set $\mbox{\boldmath $R$}=\mbox{\boldmath $W$}_e-\mbox{\boldmath $V$}$.
Since $\nabla\cdot(\mbox{\boldmath $A$}\times\mbox{\boldmath $B$})
=\nabla\times\mbox{\boldmath $A$}\cdot\mbox{\boldmath $B$}-\mbox{\boldmath $A$}\cdot\nabla\times\mbox{\boldmath $B$}$,
integration by parts yields
$$\begin{array}{c}
\displaystyle
\int_{\Bbb R^3\setminus\overline{D}}\mbox{\boldmath $R$}
\cdot\nabla\times\nabla\times\mbox{\boldmath $V$}^rdx\\
\\
\displaystyle
=-\int_{\partial D}\mbox{\boldmath $\nu$}\cdot((\nabla\times\mbox{\boldmath $V$}^r)\times\mbox{\boldmath $R$})dS
+\int_{\Bbb R^3\setminus\overline D}
\nabla\times\mbox{\boldmath $V$}^r\cdot
\nabla\times\mbox{\boldmath $R$}dx
\end{array}
\tag {4.15}
$$
and
$$\begin{array}{c}
\displaystyle
\int_{\Bbb R^3\setminus\overline D}\nabla\times\nabla\times\mbox{\boldmath $R$}\cdot\mbox{\boldmath $V$}^rdx\\
\\
\displaystyle
=-\int_{\partial D}\mbox{\boldmath $\nu$}\cdot\left((\nabla\times\mbox{\boldmath $R$})\times\mbox{\boldmath $V$}^r\right)dS
+\int_{\Bbb R^3\setminus\overline D}\nabla\times\mbox{\boldmath $R$}\cdot\nabla\times\mbox{\boldmath $V$}^rdx.
\end{array}
\tag {4.16}
$$
Taking the difference of (4.15) from (4.16) and noting $\phi\equiv 1$ in a neighbourhood of $\partial D$, we obtain
$$\begin{array}{c}
\displaystyle
\int_{\Bbb R^3\setminus\overline{D}}\left(\mbox{\boldmath $R$}
\cdot\nabla\times\nabla\times\mbox{\boldmath $V$}^r
-\nabla\times\nabla\times\mbox{\boldmath $R$}\cdot\mbox{\boldmath $V$}^r\right)dx\\
\\
\displaystyle
=\int_{\partial D}\mbox{\boldmath $\nu$}\cdot
\left((\nabla\times\mbox{\boldmath $R$})\times\mbox{\boldmath $V$}^*
-(\nabla\times\mbox{\boldmath $V$}^*)\times\mbox{\boldmath $R$}\right)dS.
\end{array}
\tag {4.17}
$$
Since $\mbox{\boldmath $R$}$ satisfies $\mbox{\boldmath $R$}\times\mbox{\boldmath $\nu$}=-\mbox{\boldmath $V$}\times\mbox{\boldmath $\nu$}$ on $\partial D$
(see (3.11)),
we have
$$\displaystyle
\mbox{\boldmath $\nu$}\cdot((\nabla\times\mbox{\boldmath $V$}^*)\times\mbox{\boldmath $R$})
=(\nabla\times\mbox{\boldmath $V$}^*)\cdot(\mbox{\boldmath $R$}\times\mbox{\boldmath $\nu$})
=-(\nabla\times\mbox{\boldmath $V$}^*)\cdot(\mbox{\boldmath $V$}\times\mbox{\boldmath $\nu$})
=-\mbox{\boldmath $V$}\cdot\mbox{\boldmath $\nu$}\times (\nabla\times\mbox{\boldmath $V$}^*).
$$
Thus, applying (4.5) to this, we obtain
$$\begin{array}{c}
\displaystyle
\mbox{\boldmath $\nu$}\cdot((\nabla\times\mbox{\boldmath $V$}^*)\times\mbox{\boldmath $R$})
=-\mbox{\boldmath $V$}\cdot\mbox{\boldmath $\nu$}\times (\nabla\times\mbox{\boldmath $V$})
=-\mbox{\boldmath $\nu$}\cdot((\nabla\times\mbox{\boldmath $V$})\times\mbox{\boldmath $V$})
\end{array}
$$
Substituting this into (2.22), we obtain
$$\begin{array}{c}
\displaystyle
J(\tau)
=\frac{1}{\mu\epsilon}\int_{\partial D}\mbox{\boldmath $\nu$}\cdot((\nabla\times\mbox{\boldmath $V$}^*)\times\mbox{\boldmath $R$})dS.
\end{array}
\tag {4.18}
$$
Moreover, from (4.4) we have
$$\begin{array}{c}
\displaystyle
\mbox{\boldmath $\nu$}\cdot((\nabla\times\mbox{\boldmath $R$})\times\mbox{\boldmath $V$}^*)
=\nabla\times\mbox{\boldmath $R$}\cdot(\mbox{\boldmath $V$}^*\times\mbox{\boldmath $\nu$})
=-\nabla\times\mbox{\boldmath $R$}\cdot(\mbox{\boldmath $V$}\times\mbox{\boldmath $\nu$})
=\mbox{\boldmath $\nu$}\times\mbox{\boldmath $V$}\cdot\nabla\times\mbox{\boldmath $R$}.
\end{array}
$$
Substituting this into (3.12), we obtain
$$\begin{array}{c}
\displaystyle
E(\tau)=
\frac{1}{\mu\epsilon}\int_{\partial D}\mbox{\boldmath $\nu$}\cdot((\nabla\times\mbox{\boldmath $R$})\times\mbox{\boldmath $V$}^*)dS
+e^{-\tau T}\int_{\Bbb R^3\setminus\overline D}\mbox{\boldmath $F$}\cdot\mbox{\boldmath $R$}dx.
\end{array}
\tag {4.19}
$$
Substituting (4.18) and (4.19) into the right-hand side on (4.17), we obtain
$$\begin{array}{c}
\displaystyle
E(\tau)-J(\tau)
=\frac{1}{\mu\epsilon}\int_{\Bbb R^3\setminus\overline{D}}\left(\mbox{\boldmath $R$}
\cdot\nabla\times\nabla\times\mbox{\boldmath $V$}^r
-\nabla\times\nabla\times\mbox{\boldmath $R$}\cdot\mbox{\boldmath $V$}^r\right)dx
+e^{-\tau T}\int_{\Bbb R^3\setminus\overline D}\mbox{\boldmath $F$}\cdot\mbox{\boldmath $R$}dx.
\end{array}
$$
From this and (3.10) we obtain
$$\begin{array}{c}
\displaystyle
E(\tau)-J(\tau)
=\int_{\Bbb R^3\setminus\overline D}\mbox{\boldmath $R$}\cdot
\left(\frac{1}{\mu\epsilon}\nabla\times\nabla\times\mbox{\boldmath $V$}^r+\tau^2\mbox{\boldmath $V$}^r\right)dx\\
\\
\displaystyle
+e^{-\tau T}\left(\int_{\Bbb R^3\setminus\overline D}\mbox{\boldmath $F$}\cdot\mbox{\boldmath $R$}\,dx
-\int_{\Bbb R^3\setminus\overline D}\mbox{\boldmath $F$}\cdot\mbox{\boldmath $V$}^rdx\right)\\
\\
\displaystyle
\equiv I+e^{-\tau T}II.
\end{array}
\tag {4.20}
$$

Here we prove
$$\displaystyle
I=O(\tau^{-1/2})J(\tau)+O(\tau^{-2}e^{-\tau T}).
\tag {4.21}
$$
From Proposition 4.2 we have
$$\begin{array}{c}
\displaystyle
\frac{1}{\mu\epsilon}\nabla\times\nabla\times\mbox{\boldmath $V$}^r(x)+\tau^2\mbox{\boldmath $V$}^r(x)
=
\phi(x)\left(\sum_{j,k,l}d_{\partial D}(x)C_{ijkl}(x)\frac{\displaystyle\partial^2\mbox{\boldmath $V$}^j}{\partial x_k\partial x_l}(x^r)\right)\\
\\
\displaystyle
+\left(\sum_{j,k,l}D_{ijkl}(x)\frac{\displaystyle\partial\mbox{\boldmath $V$}^j}{\partial x_k}(x^r)\frac{\partial\phi}{\partial x_l}(x)\right)
+\phi(x)\left(\sum_{j,k}E_{ijk}(x)\frac{\displaystyle\partial\mbox{\boldmath $V$}^j}{\partial x_k}(x^r)\right)
\\
\\
\displaystyle
+\left(\sum_{j,k,l}F_{ijkl}\mbox{\boldmath $V$}^j(x^r)\frac{\partial^2\phi}{\partial x_k\partial x_l}(x)\right),
\end{array}
$$
where $C_{ijkl}$ are of class $C^1$;
$D_{ijkl}$ and $E_{ijk}$ are of class $C^1$ and $C^0$ in a neighbourhood of $\partial D$; $F_{ijkl}$ are constants.

Substituting this into the first term on the right-hand side of (4.20) and
making a change of variables $x=y^r$, we obtain
$$\begin{array}{c}
\displaystyle
I
=\int_{D}\mbox{\boldmath $R$}(y^r)\cdot\left\{\phi(y^r)
\sum_{j,k,l}d_{\partial D}(y^r)C_{ijkl}(y^r)\frac{\displaystyle\partial^2\mbox{\boldmath $V$}^j}{\partial y_k\partial y_l}(y)
+\text{lower order terms}\right\}
J(y)dy,
\end{array}
\tag {4.22}
$$
where $J(y)$ denotes the Jacobian of the map: $y\longmapsto y^r$.
A routine involving an integration by parts and $d_{\partial D}(y^r)=d_{\partial D}(y)$ yields
$$\begin{array}{c}
\displaystyle
\int_{D}\mbox{\boldmath $R$}(y^r)\cdot\left(\phi(y^r)
\sum_{j,k,l}d_{\partial D}(y^r)C_{ijkl}(y^r)\frac{\displaystyle\partial^2\mbox{\boldmath $V$}^j}{\partial y_k\partial y_l}(y)\right)
J(y)dy
\\
\\
\displaystyle
=\left\{O(\delta)\Vert(\mbox{\boldmath $R$}^r)'\Vert_{L^2(D_{\delta})}
+O(1)\Vert\mbox{\boldmath $R$}^r\Vert_{L^2(D_{\delta})}\right\}\Vert\mbox{\boldmath $V$}'\Vert_{L^2(D_{\delta})},
\end{array}
$$
where $\mbox{\boldmath $R$}^r(y)=\mbox{\boldmath $R$}(y^r)$ and $D_{\delta}=\{y\in\,D\,\vert\,d_{\partial D}(y)<2\delta\}$.

Just simply estimating other terms on the right-hand of (4.22)
which are coming from the lower order terms, we obtain
$$\begin{array}{c}
\displaystyle
I
=O(\delta)\Vert(\mbox{\boldmath $R$}^r)'\Vert_{L^2(D_{\delta})}\Vert\mbox{\boldmath $V$}'\Vert_{L^2(D)}\\
\\
\displaystyle
+O(\delta^{-1})\Vert\mbox{\boldmath $R$}^r\Vert_{L^2(D_{\delta})}\Vert\mbox{\boldmath $V$}'\Vert_{L^2(D)}
+O(\delta^{-2})\Vert\mbox{\boldmath $R$}^r\Vert_{L^2(D_{\delta})}\Vert\mbox{\boldmath $V$}\Vert_{L^2(D)}.
\end{array}
\tag {4.23}
$$

Here we note that: using a change of variables again, we can easily obtain
$$\displaystyle
\Vert\mbox{\boldmath $R$}^r\Vert_{L^2(D_{\delta})}
\le C\Vert\mbox{\boldmath $R$}\Vert_{L^2(\Bbb R^3\setminus\overline D)},
\,\,
\Vert(\mbox{\boldmath $R$}^r)'\Vert_{L^2(D_{\delta})}
\le C\Vert\mbox{\boldmath $R$}'\Vert_{L^2((\Bbb R^3\setminus\overline D)_{\delta})},
$$
where $(\Bbb R^3\setminus\overline D)_{\delta}=\{x\in\Bbb R^3\setminus\overline D\,\vert\,
d_{\partial D}(x)<2\delta\}\subset(\Bbb R^3\setminus\overline D)_{\delta_0/2}$;
we have trivial estimates (3.26) and $\displaystyle
\Vert\mbox{\boldmath $V$}\Vert_{L^2(D)}\le \tau^{-1}\sqrt{J(\tau)}$;
from (2.15) we have $\displaystyle
\Vert\mbox{\boldmath $V$}'\Vert_{L^2(D)}\le C\sqrt{J(\tau)}$.

Therefore we see that the right-hand side on (4.23) has a bound involving $J(\tau)$, $E(\tau)$ and
$\Vert\mbox{\boldmath $R$}'\Vert_{L^2((\Bbb R^3\setminus\overline D)_{\delta})}$.

Here we describe a crucial lemma to estimate $E(\tau)$ and
$\Vert\mbox{\boldmath $R$}'\Vert_{L^2((\Bbb R^3\setminus\overline D)_{\delta})}$ in the bound.

\proclaim{\noindent Lemma 4.1.}
Assume that $\Lambda_{\partial D}(p)$ is finite and that (1.8) and (1.9) are satisfied.
Then, there exist positive constants $C$ and $\tau_0$ such that,
for all $\tau\ge\tau_0$ we have
$$\displaystyle
E(\tau)\le C(J(\tau)+e^{-2\tau T})
\tag {4.24}
$$
and
$$\displaystyle
\Vert\mbox{\boldmath $R$}'\Vert_{L^2((\Bbb R^3\setminus\overline D)_{\delta_0/2})}^2\le C(J(\tau)+e^{-2\tau T}).
\tag {4.25}
$$

\endproclaim

For the proof see the next subsection.
Once we have (4.24) and (4.25), we see that (4.23) becomes
$$\begin{array}{c}
\displaystyle
I=(O(\delta)+O((\delta\tau)^{-1})+O((\delta\tau)^{-2}))J(\tau)
+O\left((\delta\tau)^{-1})+(\delta\tau)^{-2})+\delta)\right)e^{-\tau T}\tau^{-1/2}\sqrt{J(\tau)}.
\end{array}
$$
Let $\theta>0$ and choose $\delta=\tau^{-\theta}$ with $\tau>>1$.
Then, this right-hand side becomes
$$\begin{array}{c}
\displaystyle
I
=O(\tau^{-\theta}+\tau^{-(1-\theta)}+\tau^{-2(1-\theta)})J(\tau)\\
\\
\displaystyle
+O(\tau^{-(1-\theta)}+\tau^{-2(1-\theta)}+\tau^{-\theta})e^{-\tau T}\sqrt{J(\tau)}.
\end{array}
\tag {4.26}
$$
Now choosing $\theta$ in such a way that $\theta=1-\theta$, that is, $\theta=1/2$.
Then (4.21) follows from (4.26) and (3.16).

It is easier to obtain the estimate
$\displaystyle
e^{-\tau T}II=O(e^{-\tau T}(\sqrt{E(\tau)}+\sqrt{J(\tau)}))$ than (4.21).
Then (3.16) and (3.17) give $\displaystyle e^{-\tau T}II=O(\tau^{-3/2}e^{-\tau T})$.
A combination of this and (4.21) yields
$$\displaystyle
\vert E(\tau)-J(\tau)\vert=O(\tau^{-1/2})J(\tau)+O(\tau^{-3/2}e^{-\tau T}).
\tag {4.27}
$$
Here note that (2.9) yields
$$\displaystyle
\frac{e^{-\tau T}\tau^{-3/2}}{J(\tau)}
=\frac{\displaystyle
\tau^{5+2\gamma}e^{-\tau(T-2\text{dist}\,(D,B))}\tau^{-3/2}}
{\displaystyle
\tau^{5+2\gamma}e^{2\text{dist}\,(D,B)}J(\tau)}
=O(\tau^{5+2\gamma-3/2}e^{-\tau(T-2\text{dist}\,(D,B))}).
$$
Therefore (4.27) becomes $\displaystyle \vert E(\tau)-J(\tau)\vert
=O(\tau^{-1/2})J(\tau)$. This completes the proof of (4.2).

\noindent
{\bf Remark 4.2.}
Summing up, we have shown that for the proof of (4.2) it
suffices to have estimates (4.24) and (4.25).
Note that (4.24) is sharper than (3.21).
However, we need more restrictive assumptions that $\Lambda_{\partial D}(p)$ is
finite and that (1.8) and (1.9) are satisfied.

It seems that giving an estimate of $\Vert\mbox{\boldmath $R$}'\Vert_{L^2((\Bbb
R^3\setminus\overline D)_{\delta})}$ in terms of $E(\tau)$ directly
is not trivial unlike the
scalar case.  Of course now we have (4.2) and thus, from (4.25) and (2.9) we obtain
$\Vert\mbox{\boldmath $R$}'\Vert_{L^2((\Bbb
R^3\setminus\overline D)_{\delta})}\le CE(\tau)$ for $\tau>>1$ provided $\Lambda_{\partial D}(p)$ is
finite and that (1.8) and (1.9) are satisfied.

\subsection{Proof of Lemma 4.1}

Let $\varphi$ be a smooth function on the whole space
and satisfy $0\le \varphi\le 1$; $\varphi(x)=1$ for $x$ with $d_{\partial D}(x)\le \delta_0/2$
and $\varphi(x)=0$ for $d_{\partial D}(x)\ge\delta_0$.  Here $\delta_0$ is chosen in such a way that
given $x$ with $d_{\partial D}(x)<2\delta_0$ there exists a unique $q=q(x)\in\partial D$
that attains the minimum of the function $\partial D\ni y\longmapsto \vert y-x\vert$.
We assume that $\partial D$ is $C^2$.  Then one may think that both $d_{\partial D}(x)$ and
$q(x)$ are $C^2$ for $x\in\Bbb R^3\setminus D$ with $d_{\partial D}(x)<2\delta_0$
(see \cite{GT}, p.355, Lemma 14.16).

Set $\mbox{\boldmath $R$}=\mbox{\boldmath $W$}_e-\mbox{\boldmath $V$}$.
Taking the scalar product of (3.10) with $\varphi\mbox{\boldmath $V$}^*$ and
integrating over $\Bbb R^3\setminus\overline D$, we obtain
$$\begin{array}{c}
\displaystyle
-\frac{1}{\mu\epsilon}\int_{\partial D}\mbox{\boldmath $\nu$}\times(\varphi\mbox{\boldmath $V$}^*)\cdot\nabla\times\mbox{\boldmath $R$}dS
=\frac{1}{\mu\epsilon}\int_{\Bbb R^3\setminus\overline D}\nabla\times\mbox{\boldmath $R$}\cdot\nabla\times(\varphi\mbox{\boldmath $V$}^*)dx
+\tau^2\int_{\Bbb R^3\setminus\overline D}\mbox{\boldmath $R$}\cdot(\varphi\mbox{\boldmath $V$}^*)dx\\
\\
\displaystyle
-e^{-\tau T}\int_{\Bbb R^3\setminus\overline D}\mbox{\boldmath $F$}(x,\tau)\cdot(\varphi\mbox{\boldmath $V$}^*)dx.
\end{array}
$$
Since $\mbox{\boldmath $V$}^*$ satisfies (4.4), from this and
(3.22) we obtain the expression
$$\begin{array}{c}
\displaystyle
E(\tau)
=
\frac{1}{\mu\epsilon}\int_{\Bbb R^3\setminus\overline D}\nabla\times\mbox{\boldmath $R$}\cdot\nabla\times(\varphi\mbox{\boldmath $V$}^*)dx
+\tau^2\int_{\Bbb R^3\setminus\overline D}\mbox{\boldmath $R$}\cdot(\varphi\mbox{\boldmath $V$}^*)dx\\
\\
\displaystyle
-e^{-\tau T}\left(\int_{\Bbb R^3\setminus\overline D}\mbox{\boldmath $F$}(x,\tau)\cdot(\varphi\mbox{\boldmath $V$}^*)dx
-\int_{\Bbb R^3\setminus\overline D}\mbox{\boldmath $F$}\cdot\mbox{\boldmath $R$}dx\right).
\end{array}
$$
This yields
$$\begin{array}{c}
\displaystyle
E(\tau)\le\frac{1}{\mu\epsilon}\Vert\nabla\times\mbox{\boldmath $R$}\Vert_{L^2(\Bbb R^3\setminus\overline D)}
\Vert\nabla\times(\varphi\mbox{\boldmath $V$}^*)\Vert_{L^2(\Bbb R^3\setminus\overline D)}
+\tau^2\Vert\mbox{\boldmath $R$}\Vert_{L^2(\Bbb R^3\setminus\overline D)}
\Vert\varphi\mbox{\boldmath $V$}^*\Vert_{L^2(\Bbb R^3\setminus\overline D)}\\
\\
\displaystyle
+e^{-\tau T}\Vert\mbox{\boldmath $F$}\Vert_{L^2(\Bbb R^3\setminus\overline D)}
\left(\Vert\varphi\mbox{\boldmath $V$}^*\Vert_{L^2(\Bbb R^3\setminus\overline D)}
+\Vert\mbox{\boldmath $R$}\Vert_{L^2(\Bbb R^3\setminus\overline D)}\right).
\end{array}
\tag {4.28}
$$
A change of variable $y=x-2d_{\partial D}(x)\mbox{\boldmath $\nu$}_x$ gives
$$\displaystyle
\Vert\nabla\varphi\times\mbox{\boldmath $V$}^*\Vert_{L^2(\Bbb R^3\setminus\overline D)}+
\Vert\varphi\mbox{\boldmath $V$}^*\Vert_{L^2(\Bbb R^3\setminus\overline D)}
\le C\Vert\mbox{\boldmath $V$}\Vert_{L^2(D)}\le C'\tau^{-1}\sqrt{J(\tau)}
$$
and also
$$\displaystyle
\Vert\varphi\nabla\times\mbox{\boldmath $V$}^*\Vert_{L^2(\Bbb R^3\setminus\overline D)}
\le C\Vert\mbox{\boldmath $V$}'\Vert_{L^2(D)}\le C'\sqrt{J(\tau)}.
$$
Applying these together with (3.20), (3.25) and (3.26) to the
right-hand side on (4.28), we obtain
$$\begin{array}{c}
E(\tau)
\le C\{\sqrt{E(\tau)}\sqrt{J(\tau)}+e^{-\tau T}(\sqrt{J(\tau)}+\sqrt{E(\tau)})\}.
\end{array}
$$
Now a standard argument yields (4.24).

Next we give a proof of (4.25).
Define $\mbox{\boldmath $U$}=\varphi(\mbox{\boldmath $R$}-\mbox{\boldmath $V$}^*)$.
Since $\nabla\cdot\mbox{\boldmath $R$}=e^{-\tau T}\nabla\cdot\mbox{\boldmath $F$}/\tau^2$ and $\nabla\cdot
\mbox{\boldmath $F$}=-\tau\nabla\cdot\mbox{\boldmath $E$}(x,T)$, we have
$\displaystyle
\nabla\cdot\mbox{\boldmath $R$}=-e^{-\tau T}\tau^{-1}\nabla\cdot\mbox{\boldmath $E$}(x,T)$.
However, from the governing equations of $\mbox{\boldmath $E$}$ and $\mbox{\boldmath $H$}$ in the time domain
we have
$$
\displaystyle
\nabla\cdot\mbox{\boldmath $E$}(x,T)
=\int_0^T\nabla\cdot\mbox{\boldmath $J$}(x,t)dt.
$$
So choosing $\delta_0$ in such a way that
$\overline B\cap(\Bbb R^3\setminus\overline D)_{\delta_0/2}=\emptyset$,
we conclude $\nabla\cdot\mbox{\boldmath $R$}=0$ in $(\Bbb R^3\setminus\overline D)_{\delta_0/2}$
and hence
$\nabla\cdot\mbox{\boldmath $U$}=\nabla\varphi\cdot\mbox{\boldmath $R$}-\nabla\cdot(\varphi\mbox{\boldmath $V$}^*)
\in L^2((\Bbb R^3\setminus\overline D)_{\delta_0/2})$.
Moreover, $\mbox{\boldmath $U$}$ together with $\nabla\times\mbox{\boldmath $U$}$ belongs to $L^2((\Bbb R^3\setminus\overline D)_{\delta_0/2})$;
$\mbox{\boldmath $U$}$ satisfies $\mbox{\boldmath $U$}\times\mbox{\boldmath $\nu$}=\mbox{\boldmath $0$}$ on $\partial D$
and $\mbox{\boldmath $U$}=\mbox{\boldmath $0$}$ in $(\Bbb R^3\setminus\overline D)_{\delta_0/2}$.
Therefore, by Corollary 1.1 on p. 212 and (ii) of Remark 2 on p. 213 in \cite{DL3},
we have $U\in H^1((\Bbb R^3\setminus\overline D)_{\delta_0/2})$
and
$$\displaystyle
\Vert\mbox{\boldmath $U$}\Vert_{H^1((\Bbb R^3\setminus\overline D)_{\delta_0/2})}^2
\le C(\Vert\mbox{\boldmath $U$}\Vert_{L^2((\Bbb R^3\setminus\overline D)_{\delta_0/2})}^2
+\Vert\nabla\times\mbox{\boldmath $U$}\Vert_{L^2((\Bbb R^3\setminus\overline D)_{\delta_0/2})}^2
+\Vert\nabla\cdot\mbox{\boldmath $U$}\Vert_{L^2((\Bbb R^3\setminus\overline D)_{\delta_0/2})}^2).
$$
Applying (3.25) and (3.26) and a change of variables to this right-hand side, we obtain
$$\displaystyle
\Vert\mbox{\boldmath $U$}\Vert_{H^1(\Bbb R^3\setminus\overline D)}^2
\le C(E(\tau)+
\Vert\mbox{\boldmath $V$}\Vert_{H^1(D)}^2).
\tag {4.29}
$$
Since $\varphi\mbox{\boldmath $R$}=\mbox{\boldmath $U$}+\varphi\mbox{\boldmath $V$}^*$,
(4.29) together with the estimate $\Vert\varphi\mbox{\boldmath $V$}^*\Vert_{H^1(\Bbb R^3\setminus\overline D)}
\le C\Vert\mbox{\boldmath $V$}\Vert_{H^1(D)}\le C'\sqrt{J(\tau)}$ gives
$\displaystyle
\Vert\mbox{\boldmath $R$}'\Vert_{L^2((\Bbb R^3\setminus\overline D)_{\delta_0/2})}^2
\le C(E(\tau)+J(\tau))$.  A combination of this and (4.24) yields (4.25).
This completes the proof of Lemma 4.1.

\section{Concluding remarks}

In this paper, we employed a simple form (1.2) as a model of the current density, however, in principle,
it may be possible to cover more complicated model of the current density, at least, in the frame work
of the solution as constructed in \cite{DL5}.

The method presented here can be applied also to an {\it interior problem} similar
to that considered in \cite{IE3}.  The problem therein aims at extracting information about the geometry
of an unknown cavity from the wave in time domain which is produced by the initial data
localized inside the cavity and propagates therein.

Single measurement version of the time domain enclosure method also finds an application
to an inverse initial boundary value problem for the heat equation in three-space
dimensions.  For this see Theorem 1.1 in \cite{IK} and consult Section 3 in \cite{IH} for an open problem in the visco elasticity.

Some of open problems are in order.

$\bullet$  A lot of papers deals with the perfectly conducting obstacle as the first step
(see \cite{CK} and references therein). It is a typical and important condition as everyone first considers,
like the Dirichlet boundary condition for the wave equation.
This paper also follows that traditional order
and note that the aim of this paper is to introduce a method for inverse electromagnetic obstacle scattering.
However, as a next step, it is natural to ask:
how about the case when the electromagnetic wave satisfies a more general boundary condition
like the Leontovich condition on the surface of the obstacle(see, e.g., \cite{ALN})?
Note that, for the wave equation with the Robin type boundary condition
we have \cite{IE2} and \cite{IE5} which contain results corresponding to Theorem 1.1 and Theorem 1.2,
respectively.

$\bullet$  How about the case when the reflected electromagnetic wave is observed at a different place
from the support of the source ?
We expect that the observed data give us different information about the geometry of unknown obstacle
together with a constructive method which yields the location of all the first reflection points
from a single observed wave as seen for an acoustic wave case in \cite{IE4}.
It would be interesting to see also \cite{BGJ} for a comparison of monostatic and bistatic radar images.

$\bullet$  There are several other inverse obstacle scattering
problems in time domain whose governing equations are {\it systems
of partial differential equations}.  Extend the range of the
applications of the method presented here to such systems.

$$\quad$$

\centerline{{\bf Acknowledgement}}

This research was partially supported by Grant-in-Aid for
Scientific Research (C)(No. 25400155) of Japan  Society for the
Promotion of Science.

\section{Appendix}

\subsection{Proof of Claim 1.}

Clearly $\tilde{\mbox{\boldmath $V$}}$ satisfies (4.8).
To check (4.9) we have to compute $\nabla\times\tilde{\mbox{\boldmath $V$}}$.
Set $\displaystyle B(x)=\mbox{\boldmath $V$}(x)\cdot\mbox{\boldmath $n$}(x)$.

We have
$$\displaystyle
\nabla(B(x^r))
=(2q'(x)^T-I)(\nabla B)(x^r).
\tag {A.1}
$$
Let $y\in\partial D$.
Since $q'(y)\mbox{\boldmath $\nu$}_{y}=0$, we get
$$\displaystyle
\nabla(B(x^r))\vert_{x=y}\cdot\mbox{\boldmath $\nu$}_{y}=-(\nabla B)(y)\cdot\mbox{\boldmath $\nu$}_{y}.
$$
On the other hand, we have $q'(y)\mbox{\boldmath $v$}=\mbox{\boldmath $v$}$ for all vectors with
$\mbox{\boldmath $\nu$}_{y}\cdot\mbox{\boldmath $v$}=0$.
Thus (A.1) gives
$$\displaystyle
\nabla(B(x^r))\vert_{x=y}\cdot\mbox{\boldmath $v$}=(\nabla B)(y)\cdot\mbox{\boldmath $v$}.
$$
From these we obtain
$$\displaystyle
\nabla(B(x^r))\vert_{x=y}=\{(\nabla B)(y)-((\nabla B)(y)\cdot\mbox{\boldmath $\nu$}_{y})
\mbox{\boldmath $\nu$}_{y}\}-((\nabla B)(y)\cdot\mbox{\boldmath $\nu$}_{y})\mbox{\boldmath $\nu$}_{y}.
\tag {A.2}
$$
Here we note that $\mbox{\boldmath $n$}(x)=\nabla(d_{\partial D}(x))$ for $x\in\Bbb R^3\setminus\overline D$
and $\mbox{\boldmath $n$}(x)=-\nabla (d_{\partial D}(x))$ for $x\in D$.
This gives $\nabla\times\mbox{\boldmath $n$}=\mbox{\boldmath $0$}$.  Thus, we have
$\displaystyle
(\nabla\times\tilde{\mbox{\boldmath $B$}})(x)
=\nabla(B(x^r))\times\mbox{\boldmath $n$}(x)$
and $\displaystyle
(\nabla\times\mbox{\boldmath $B$})(x)
=(\nabla B)(x)\times\mbox{\boldmath $n$}(x)$,
where $\tilde{\mbox{\boldmath $B$}}(x)=\mbox{\boldmath $B$}(x^r)$.
Thus, from (A.2) we obtain
$$\displaystyle
\nabla\times\tilde{\mbox{\boldmath $B$}}
=\nabla\times\mbox{\boldmath $B$}\,\,\text{on}\,\partial D.
\tag {A.3}
$$

Define $\tilde{\mbox{\boldmath $A$}}(x)=-\mbox{\boldmath $A$}(x^r)$ for $x\in\Bbb R^3\setminus\overline D$.
Let $y\in\partial D$.
Applying (A.2) for $B$ replaced with $-\mbox{\boldmath $A$}^i$ for each $i=1,2,3$, we have
$$\begin{array}{c}
\displaystyle
(\nabla\tilde{\mbox{\boldmath $A$}}^i)(y)
=-(\nabla\mbox{\boldmath $A$}^i)(y)+2((\nabla\mbox{\boldmath $A$}^i)(y)\cdot\mbox{\boldmath $\nu$}_{y})
\mbox{\boldmath $\nu$}_{y}.
\end{array}
$$
Note that
$$\begin{array}{c}
\displaystyle
\nabla\times\tilde{\mbox{\boldmath $A$}}
=\sum_{i=1}^3\nabla\times(\tilde{\mbox{\boldmath $A$}}^i\mbox{\boldmath $e$}_i)
=\sum_{i=1}^3\nabla\tilde{\mbox{\boldmath $A$}}^i\times\mbox{\boldmath $e$}_i
\end{array}
$$
and the same for $\nabla\times\mbox{\boldmath $A$}$.  These yield
$$\begin{array}{c}
\displaystyle
(\nabla\times\tilde{\mbox{\boldmath $A$}})(y)
=-(\nabla\times\mbox{\boldmath $A$})(y)
+2\sum_{i=1}^3\{(\nabla\mbox{\boldmath $A$}^i)(y)\cdot\nu_{y}\}\mbox{\boldmath $\nu$}_{y}\times\mbox{\boldmath $e$}_i.
\end{array}
\tag {A.4}
$$
Write
$$\begin{array}{c}
\displaystyle
\sum_{i=1}^3\{(\nabla\mbox{\boldmath $A$}^i)(y)\cdot\nu_{y}\}\mbox{\boldmath $\nu$}_{y}\times\mbox{\boldmath $e$}_i
=\mbox{\boldmath $\nu$}_{y}\times\sum_{i=1}^3\{(\nabla\mbox{\boldmath $A$}^i(y)\cdot\mbox{\boldmath $\nu$}_{y}\}\mbox{\boldmath $e$}_i
=\mbox{\boldmath $\nu$}_{y}\times\{\mbox{\boldmath $A$}'(y)
\mbox{\boldmath $\nu$}_{y}\}.
\end{array}
$$
Then, (A.4) becomes
$$\displaystyle
(\nabla\times\tilde{\mbox{\boldmath $A$}})(y)
=-(\nabla\times\mbox{\boldmath $A$})(y)
+2\mbox{\boldmath $\nu$}_y\times\{\mbox{\boldmath $A$}'(y)\mbox{\boldmath $\nu$}_{y}\}.
\tag {A.5}
$$
Taking the vector product of both sides on (A.5) with $\mbox{\boldmath $\nu$}_{y}$, we obtain
$$\begin{array}{c}
\displaystyle
\mbox{\boldmath $\nu$}_{y}\times(\nabla\times\tilde{\mbox{\boldmath $A$}})(y)
=-\mbox{\boldmath $\nu$}_{y}\times(\nabla\times\mbox{\boldmath $A$})(y)-2\mbox{\boldmath $A$}'(y)
\mbox{\boldmath $\nu$}_{y}.
\end{array}
\tag {A.6}
$$
Note that, in the derivation of this,  we have made use of
the identity
$$\displaystyle
\mbox{\boldmath $A$}\times(\mbox{\boldmath $B$}\times\mbox{\boldmath $C$})
=\mbox{\boldmath $B$}(\mbox{\boldmath $A$}\cdot\mbox{\boldmath $C$})-\mbox{\boldmath $C$}(\mbox{\boldmath $A$}\cdot
\mbox{\boldmath $B$});
\tag {A.7}
$$
the equation $\mbox{\boldmath $A$}'(x^r)\mbox{\boldmath $n$}(x)\cdot\mbox{\boldmath $n$}(x)=0$
which is an easy consequence of the property $\mbox{\boldmath $A$}(x)\cdot\mbox{\boldmath $n$}(x)=0$.

It is easy to see that, for any vector $\mbox{\boldmath $v$}$, we have
$$\displaystyle
\mbox{\boldmath $v$}\times(\nabla\times\mbox{\boldmath $A$})(x)
=(\mbox{\boldmath $A$}'(x)^T-\mbox{\boldmath $A$}'(x))\mbox{\boldmath $v$}.
\tag {A.8}
$$
Rewrite the right-hand side on (A.6) as
$$\displaystyle
\mbox{\boldmath $\nu$}_{y}\times(\nabla\times\mbox{\boldmath $A$})(y)
-2\left\{\mbox{\boldmath $\nu$}_{y}\times(\nabla\times\mbox{\boldmath $A$})(y)+\mbox{\boldmath $A$}'(y)
\mbox{\boldmath $\nu$}_{y}\right\}.
$$
Then applying (A.8) to the second term of this, we know that (A.6) becomes
$$\displaystyle
\mbox{\boldmath $\nu$}_{y}\times(\nabla\times\tilde{\mbox{\boldmath $A$}})(y)
=
\mbox{\boldmath $\nu$}_{y}\times(\nabla\times\mbox{\boldmath $A$})(y)
-2(\mbox{\boldmath $A$}'(y))^T
\mbox{\boldmath $\nu$}_{y}.
\tag {A.9}
$$
Let $y(\sigma)$ be an arbitrary curve on $\partial D$ with
$y(0)=y$. We have $\mbox{\boldmath
$A$}(y(\sigma))\cdot\mbox{\boldmath $\nu$}_{y(\sigma)}=0$. Differentiating this both
sides with respect to $\sigma_j$, we obtain,
$$\displaystyle
\mbox{\boldmath $A$}'(y)\frac{\partial y}{\partial\sigma_j}\vert_{\sigma=0}\cdot
\mbox{\boldmath $\nu$}_y
=-\mbox{\boldmath $A$}(y)\cdot\frac{\partial}{\partial\sigma_j}(\mbox{\boldmath $\nu$}_{y(\sigma)})\vert_{\sigma=0}.
$$
Recalling the definition and symmetry of the shape operator for $\partial D$ at $y\in\partial D$ with respect to
$\mbox{\boldmath $\nu$}_y$,
we have
$$\displaystyle
\mbox{\boldmath $A$}'(y)^T\mbox{\boldmath $\nu$}_y\cdot\mbox{\boldmath $v$}
=S_y(\partial D)\mbox{\boldmath $A$}(y)\cdot\mbox{\boldmath $v$}
\tag {A.10}
$$
for all tangent vectors $\mbox{\boldmath $v$}$
at $y$ of $\partial D$.  Since $S_y(\partial D)\mbox{\boldmath $A$}(x)$ is a tangent vector at $y$ of $\partial D$
and $\mbox{\boldmath $A$}'(y)^T\mbox{\boldmath $\nu$}_y\cdot\mbox{\boldmath $\nu$}_y=0$,
we know that (A.10) is valid for all vectors of $\Bbb R^3$.  Thus we obtain
$\displaystyle
\mbox{\boldmath $A$}'(y)^T\mbox{\boldmath $\nu$}_y=S_y(\partial D)\mbox{\boldmath $A$}(y)$.
Now from this and (A.9) we obtain
$$\displaystyle
\mbox{\boldmath $\nu$}_x\times(\nabla\times\tilde{\mbox{\boldmath $A$}})(x)
=\mbox{\boldmath $\nu$}_x\times(\nabla\times\mbox{\boldmath $A$})(x)
-2S_x(\partial D)(\mbox{\boldmath $A$}(x))\,\,\text{on}\,\partial D.
$$
Now from this and (A.3) we obtain (4.9).

\noindent
$\Box$

\subsection{Proof of Claim 2.}

Since $\nabla(d_{\partial D}(x))=\mbox{\boldmath $n$}(x)$,
we have $\mbox{\boldmath $n$}'=(\mbox{\boldmath $n$}')^T$
and thus $\mbox{\boldmath $n$}'(x)\mbox{\boldmath $n$}(x)=\mbox{\boldmath $0$}$.
These yield $\displaystyle
\mbox{\boldmath $C$}(x)=2d_{\partial D}(x)\mbox{\boldmath $n$}'\mbox{\boldmath $A$}(x^r)$
and $\displaystyle
\nabla\times\mbox{\boldmath $C$}(x)
=2\mbox{\boldmath $n$}\times (\mbox{\boldmath $n$}'\mbox{\boldmath $A$}(x^r))
+2d_{\partial D}(x)\nabla\times(\mbox{\boldmath $n$}'\mbox{\boldmath $A$}(x^r))$.
Thus
$\displaystyle
\nabla\times\mbox{\boldmath $C$}
=2\mbox{\boldmath $n$}\times (\mbox{\boldmath $n$}'\mbox{\boldmath $A$})\,\,\text{on}\,\partial D$.
Using (A.7) and $(\mbox{\boldmath $n$}')^T\mbox{\boldmath $n$}=\mbox{\boldmath $0$}$, from this we obtain
$$\begin{array}{c}
\displaystyle
\mbox{\boldmath $n$}\times(\nabla\times\mbox{\boldmath $C$})
=2\mbox{\boldmath $n$}\times\{\mbox{\boldmath $n$}\times (\mbox{\boldmath $n$}'\mbox{\boldmath $A$})\}
=2\mbox{\boldmath $n$}(\mbox{\boldmath $n$}\cdot \mbox{\boldmath $n$}'\mbox{\boldmath $A$})-
2\mbox{\boldmath $n$}'\mbox{\boldmath $A$}(\mbox{\boldmath $n$}\cdot\mbox{\boldmath $n$})
=-2\mbox{\boldmath $n$}'\mbox{\boldmath $A$}.
\end{array}
$$
Since $-\mbox{\boldmath $n$}'\mbox{\boldmath $A$}=S(\partial D)\mbox{\boldmath $A$}$ on $\partial D$, we obtain (4.10).

\noindent
$\Box$

\subsection{Proof of Proposition 4.2.}

It is clear that Proposition 4.2 is a direct consequence of (A.20) and (A.22) in the following subsubsections.

\subsubsection{Computation of $(1/\mu\epsilon)\nabla\times\nabla\times\tilde{\mbox{\boldmath $V$}}+\tau^2\tilde{\mbox{\boldmath $V$}}$
for $\tilde{\mbox{\boldmath $V$}}$ given by (4.6).}

We have
$$\begin{array}{c}
\displaystyle
(\nabla\cdot\tilde{\mbox{\boldmath $V$}})(x)
=-\nabla\cdot(\mbox{\boldmath $A$}(x^r))+\nabla\cdot(\mbox{\boldmath $B$}(x^r))\\
\\
\displaystyle
=-(\nabla\cdot\mbox{\boldmath $A$})(x^r)-2d_{\partial D}(x)\text{Trace}\,(\mbox{\boldmath $A$}'(x^r)\mbox{\boldmath $n$}'(x))\\
\\
\displaystyle
+(\nabla\cdot\mbox{\boldmath $B$})(x^r)-
2\mbox{\boldmath $B$}'(x^r)\mbox{\boldmath $n$}(x)\cdot\mbox{\boldmath $n$}(x)
-2d_{\partial D}(x)\text{Trace}\,(\mbox{\boldmath $B$}'(x^r)\mbox{\boldmath $n$}'(x))\\
\\
\displaystyle
=-(\nabla\cdot\mbox{\boldmath $V$})(x^r)
+2\left((\nabla\cdot\mbox{\boldmath $B$})(x^r)
-\mbox{\boldmath $B$}'(x^r)\mbox{\boldmath $n$}(x)\cdot\mbox{\boldmath $n$}(x)\right)\\
\\
\displaystyle
-2d_{\partial D}(x)\text{Trace}\,(\mbox{\boldmath $V$}'(x^r)\mbox{\boldmath $n$}'(x)).
\end{array}
$$
Since $\mbox{\boldmath $n$}'\mbox{\boldmath $n$}=\mbox{\boldmath $0$}$ and $\displaystyle
\mbox{\boldmath $B$}'=B\mbox{\boldmath $n$}'+\mbox{\boldmath $n$}\otimes\nabla B$,
we have $\mbox{\boldmath $B$}'\mbox{\boldmath $n$}=(\nabla B\cdot\mbox{\boldmath $n$})
\mbox{\boldmath $n$}$ and thus $\mbox{\boldmath $B$}'(x^r)\mbox{\boldmath $n$}(x)\cdot
\mbox{\boldmath $n$}(x)=\nabla B(x^r)\cdot\mbox{\boldmath $n$}(x)$.
Since $\displaystyle
\nabla\cdot\mbox{\boldmath $B$}
=\nabla B\cdot\mbox{\boldmath $n$}+B\nabla\cdot\mbox{\boldmath $n$}$, we obtain
$$
\displaystyle
(\nabla\cdot\mbox{\boldmath $B$})(x^r)
-\mbox{\boldmath $B$}'(x^r)\mbox{\boldmath $n$}(x)\cdot\mbox{\boldmath $n$}(x)
=B(x^r)(\nabla\cdot\mbox{\boldmath $n$})(x^r).
$$
Since $\mbox{\boldmath $V$}$ satisfies (4.11), taking the rotation of the both sides,
one gets
$$\displaystyle
\nabla\cdot\mbox{\boldmath $V$}=0\,\,\text{in}\,D.
\tag {A.11}
$$
From these we obtain
$$\begin{array}{c}
\displaystyle
(\nabla\cdot\tilde{\mbox{\boldmath $V$}})(x)
=2B(x^r)(\nabla\cdot\mbox{\boldmath $n$})(x^r)
-2d_{\partial D}(x)
\text{Trace}\,\left(\mbox{\boldmath $V$}'(x^r)\mbox{\boldmath $n$}'(x)\right).
\end{array}
$$
Further a direct computation yields
$$\begin{array}{c}
\displaystyle
\nabla\{B(x^r)(\nabla\cdot\mbox{\boldmath $n$})(x^r)\}
=(\nabla\cdot\mbox{\boldmath $n$})(x^r)(I-2\pi(x))(\nabla B)(x^r)\\
\\
\displaystyle
+B(x^r)
(I-2\pi(x))(\nabla(\nabla\cdot\mbox{\boldmath $n$}))(x^r)-2d_{\partial D}(x)\mbox{\boldmath $n$}'(x)
\{\nabla(B\nabla\cdot\mbox{\boldmath $n$})\}(x^r)
\end{array}
$$
and
$$\begin{array}{c}
\displaystyle
\nabla\{\text{Trace}\,\left(\mbox{\boldmath $V$}'(x^r)\mbox{\boldmath $n$}'(x)\right)\}
=R^{2,0}(x)\nabla^2\mbox{\boldmath $V$}(x^r)+R^{1,0}(x)\nabla\mbox{\boldmath $V$}(x^r)
-2d_{\partial D}(x)R^{2,1}(x)\nabla^2\mbox{\boldmath $V$}(x^r),
\end{array}
$$
where
$$\begin{array}{c}
\displaystyle
R^{2,0}(x)\nabla^2\mbox{\boldmath $V$}(x^r)
=\left(\sum_{i,k,l}(\delta_{lj}-2n_ln_j)\frac{\partial n^k}{\partial x_i}(x)\frac{\partial^2\mbox{\boldmath $V$}^i}
{\partial x_l\partial x_k}(x^r)\right),
\end{array}
$$
$$\begin{array}{c}
\displaystyle
R^{2,1}(x)\nabla^2\mbox{\boldmath $V$}(x^r)
=\left(\sum_{i,k,l}\frac{\partial n^l}{\partial x_j}(x)\frac{\partial n^k}{\partial x_i}(x)\frac{\partial^2\mbox{\boldmath $V$}^i}
{\partial x_l\partial x_k}(x^r)\right)
\end{array}
$$
and
$$\begin{array}{c}
\displaystyle
R^{1,0}\nabla\mbox{\boldmath $V$}(x^r)
=\left(\sum_{i,k}\frac{\partial^2 n^k}{\partial x_j\partial x_i}(x)
\frac{\partial\mbox{\boldmath $V$}^i}
{\partial x_k}(x^r)\right).
\end{array}
$$
From these we obtain
$$\begin{array}{c}
\displaystyle
\nabla(\nabla\cdot\tilde{\mbox{\boldmath $V$}})(x)
\\
\\
\displaystyle
=2(\nabla\cdot\mbox{\boldmath $n$})(x^r)(I-2\pi(x))(\nabla B)(x^r)
+2B(x^r)
(I-2\pi(x))(\nabla(\nabla\cdot\mbox{\boldmath $n$}))(x^r)\\
\\
\displaystyle
-2\mbox{\boldmath $n$}(x)\text{Trace}\,\left(\mbox{\boldmath $V$}'(x^r)\mbox{\boldmath $n$}'(x)\right)\}
-2d_{\partial D}(x)\mbox{\boldmath $Z$}(x,\nabla^2\mbox{\boldmath $V$}(x^r),\nabla\mbox{\boldmath $V$}(x^r),\mbox{\boldmath $V$}(x^r)),
\end{array}
\tag {A.12}
$$
where
$$\begin{array}{c}
\displaystyle
\mbox{\boldmath $Z$}(x,\nabla^2\mbox{\boldmath $V$}(x^r),\nabla\mbox{\boldmath $V$}(x^r),\mbox{\boldmath $V$}(x^r))\\
\\
\displaystyle
=2\mbox{\boldmath $n$}'(x)\{\nabla(B\nabla\cdot\mbox{\boldmath $n$})\}(x^r)
+R^{2,0}(x)\nabla^2\mbox{\boldmath $V$}(x^r)+R^{1,0}(x)\nabla\mbox{\boldmath $V$}(x^r)\\
\\
\displaystyle
-2d_{\partial D}(x)R^{2,1}(x)\nabla^2\mbox{\boldmath $V$}(x^r).
\end{array}
$$

A direct computation yields
$$\begin{array}{c}
\displaystyle
\triangle(\mbox{\boldmath $B$}(x^r))
=(\triangle\mbox{\boldmath $B$})(x^r)-2(\nabla\cdot\mbox{\boldmath $n$})(x)((\nabla B)(x^r)
\cdot\mbox{\boldmath $n$}(x))\mbox{\boldmath $n$}(x)\\
\\
\displaystyle
+2(\mbox{\boldmath $n$}'(x)-\mbox{\boldmath $n$}'(x^r))(\nabla B)(x^r)+B(x^r)((\triangle\mbox{\boldmath $n$})(x)
-(\triangle\mbox{\boldmath $n$})(x^r))
\\
\\
\displaystyle
-2d_{\partial D}(x)\mbox{\boldmath $N$}(x,\nabla^2B(x^r),\nabla B(x^r)),
\end{array}
\tag {A.13}
$$
where
$$\begin{array}{c}
\displaystyle
\mbox{\boldmath $N$}(x,\nabla^2B(x^r),\nabla B(x^r))\\
\\
\displaystyle
=
\left(\text{Trace}\,\{(\nabla^2 B)(x^r)\mbox{\boldmath $n$}'(x)\}
-2d_{\partial D}(x)\text{Trace}\{(\mbox{\boldmath $n$}'(x))^2(\nabla^2 B)(x^r)\}\right)
\mbox{\boldmath $n$}(x)
\\
\\
\displaystyle
+\{\nabla(\nabla\cdot\mbox{\boldmath $n$})(x)\cdot(\nabla B)(x^r)\}
\mbox{\boldmath $n$}(x)+2\mbox{\boldmath $n$}'(x)^2(\nabla B)(x^r).
\end{array}
$$
And also
$$\begin{array}{c}
\displaystyle
\triangle(\mbox{\boldmath $A$}(x^r))
=(\triangle\mbox{\boldmath $A$})(x^r)-2(\nabla\cdot\mbox{\boldmath $n$})(x)
\mbox{\boldmath $n$}'(x)\mbox{\boldmath $A$}(x^r)\\
\\
\displaystyle
-2d_{\partial D}(x)\mbox{\boldmath $T$}(x,(\nabla^2\mbox{\boldmath $A$})(x^r),(\nabla\mbox{\boldmath $A$})(x^r)),
\end{array}
\tag {A.14}
$$
where the $i$-th component of $\mbox{\boldmath $T$}(x,\nabla^2\mbox{\boldmath $A$}(x^r),\nabla\mbox{\boldmath $A$}(x^r))$
is given by
$$\begin{array}{c}
\displaystyle
\left(\mbox{\boldmath $T$}(x,(\nabla^2\mbox{\boldmath $A$})(x^r),(\nabla\mbox{\boldmath $A$})(x^r))\right)^i
\\
\\
\displaystyle
=\text{Trace}\,\{(\nabla\mbox{\boldmath $A$}^i)'(x^r)\mbox{\boldmath $n$}'(x)\}
+\text{Trace}\,\left\{\mbox{\boldmath $n$}'(x)(\nabla^2\mbox{\boldmath $A$}^i)(x^r)\right\}\\
\\
\displaystyle
-d_{\partial D}(x)
\text{Trace}\,\left\{(\mbox{\boldmath $n$}'(x))^2(\nabla^2\mbox{\boldmath $A$}^i)(x^r)\right\}
+\left(\mbox{\boldmath $A$}'(x^r)\{\nabla(\nabla\cdot\mbox{\boldmath $n$})\}(x)\right)^i.
\end{array}
$$

Using (4.11), (A.11) and the formula
$\displaystyle
\nabla\times\nabla\times\mbox{\boldmath $V$}=
\nabla(\nabla\cdot\mbox{\boldmath $V$})-\triangle\mbox{\boldmath $V$}$,
we have
$$\displaystyle
-\frac{1}{\mu\epsilon}\triangle\mbox{\boldmath $V$}+\tau^2\mbox{\boldmath $V$}=\mbox{\boldmath $0$}\,\,\text{in}\, D.
\tag {A.15}
$$
From this together with (4.13) and (4.6) we have
$$\begin{array}{c}
\displaystyle
-\frac{1}{\mu\epsilon}\triangle\tilde{\mbox{\boldmath $V$}}
+\tau^2\tilde{\mbox{\boldmath $V$}}
=-\frac{1}{\mu\epsilon}(\triangle\tilde{\mbox{\boldmath $V$}}(x)+(\triangle\mbox{\boldmath $V$})(x^r))
+\tau^2(\triangle\tilde{\mbox{\boldmath $V$}}(x)+(\triangle\mbox{\boldmath $V$})(x^r))\\
\\
\displaystyle
=2\left(-\frac{1}{\mu\epsilon}(\triangle\mbox{\boldmath $B$})(x^r)+\tau^2\mbox{\boldmath $B$}(x^r)\right)+\mbox{\boldmath $R$}(x),
\end{array}
\tag {A.16}
$$
where
$$\begin{array}{c}
\displaystyle
\mbox{\boldmath $R$}(x)
=\frac{1}{\mu\epsilon}\left\{(\triangle(\mbox{\boldmath $A$}(x^r))-(\triangle\mbox{\boldmath $A$})(x^r))
-
(\triangle(\mbox{\boldmath $B$}(x^r))
-(\triangle\mbox{\boldmath $B$})(x^r))\right\}.
\end{array}
\tag {A.17}
$$
Here we claim the expression
$$\begin{array}{c}
\displaystyle
-\frac{1}{\mu\epsilon}(\triangle\mbox{\boldmath $B$})(x^r)+\tau^2\mbox{\boldmath $B$}(x^r)\\
\\
\displaystyle
=-\frac{1}{\mu\epsilon}\left\{\mbox{\boldmath $n$}'(x^r)(\nabla B)(x^r)+B(x^r)(I-\pi(x))(\triangle\mbox{\boldmath $n$})(x^r)\right\}\\
\\
\displaystyle
-\frac{1}{\mu\epsilon}\left(2\text{Trace}\,(\mbox{\boldmath $A$}'(x^r)\mbox{\boldmath $n$}'(x^r))
+\mbox{\boldmath $A$}(x^r)\cdot(\triangle\mbox{\boldmath $n$})(x^r)\right)\mbox{\boldmath $n$}(x).
\end{array}
\tag {A.18}
$$
This is proved as follows.

From (4.13) and (A.15), we have
$$\displaystyle
-\frac{1}{\mu\epsilon}\triangle\mbox{\boldmath $B$}+\tau^2\mbox{\boldmath $B$}
=\frac{1}{\mu\epsilon}\triangle\mbox{\boldmath $A$}-\tau^2\mbox{\boldmath $A$}.
\tag {A.19}
$$
This gives
$$\displaystyle
\left(-\frac{1}{\mu\epsilon}(\triangle\mbox{\boldmath $B$})(x^r)+\tau^2\mbox{\boldmath $B$}(x^r)\right)\cdot\mbox{\boldmath $n$}(x)
=\frac{1}{\mu\epsilon}(\triangle\mbox{\boldmath $A$})(x^r)\cdot\mbox{\boldmath $n$}(x).
$$
Since $\mbox{\boldmath $A$}(x)\cdot\mbox{\boldmath $n$}(x)=0$, we have
$$\displaystyle
\triangle\mbox{\boldmath $A$}(x)\cdot\mbox{\boldmath $n$}(x)
=-2\text{Trace}\,(\mbox{\boldmath $A$}'(x)(\mbox{\boldmath $n$}'(x))^T)-\mbox{\boldmath $A$}(x)\cdot(\triangle\mbox{\boldmath $n$})(x).
$$
Since $\mbox{\boldmath $n$}(x^r)=\mbox{\boldmath $n$}(x)$ and $\mbox{\boldmath $n$}'(x)^T=\mbox{\boldmath $n$}'(x)$, from this we obtain
$$\displaystyle
(\triangle\mbox{\boldmath $A$})(x^r)\cdot\mbox{\boldmath $n$}(x)
=-2\text{Trace}\,(\mbox{\boldmath $A$}'(x^r)\mbox{\boldmath $n$}'(x^r))-\mbox{\boldmath $A$}(x^r)\cdot(\triangle\mbox{\boldmath $n$})(x^r).
$$
Thus, we have
$$\displaystyle
\pi(x)\left(-\frac{1}{\mu\epsilon}(\triangle\mbox{\boldmath $B$})(x^r)+\tau^2\mbox{\boldmath $B$}(x^r)\right)
=-\frac{1}{\mu\epsilon}\left(2\text{Trace}\,(\mbox{\boldmath $A$}'(x^r)\mbox{\boldmath $n$}'(x^r))
+\mbox{\boldmath $A$}(x^r)\cdot(\triangle\mbox{\boldmath $n$})(x^r)\right)
\mbox{\boldmath $n$}(x).
$$
On the other hand,
$$\displaystyle
(I-\pi(x))\left(-\frac{1}{\mu\epsilon}(\triangle\mbox{\boldmath $B$})(x^r)+\tau^2\mbox{\boldmath $B$}(x^r)\right)
=-\frac{1}{\mu\epsilon}(I-\pi(x))(\triangle\mbox{\boldmath $B$})(x^r).
$$
Since
$$\displaystyle
\triangle\mbox{\boldmath $B$}(x)
=(\triangle B)(x)\mbox{\boldmath $n$}(x)+\mbox{\boldmath $n$}'(x)(\nabla B)(x)+B(x)(\triangle\mbox{\boldmath $n$})(x),
$$
and $\mbox{\boldmath $n$}'(x^r)^T\mbox{\boldmath $n$}(x)=\mbox{\boldmath $0$}$,
we obtain
$$\begin{array}{c}
\displaystyle
(I-\pi(x))\left(-\frac{1}{\mu\epsilon}(\triangle\mbox{\boldmath $B$})(x^r)+\tau^2\mbox{\boldmath $B$}(x^r)\right)\\
\\
\displaystyle
=-\frac{1}{\mu\epsilon}\left\{\mbox{\boldmath $n$}'(x^r)(\nabla B)(x^r)+B(x^r)(I-\pi(x))(\triangle\mbox{\boldmath $n$})(x^r)\right\}.
\end{array}
$$
Summing up, we obtain (A.18).

Now from (A.12), (A.13), (A.14), (A.16), (A.17), (A.18)
and using the formula $\nabla\times\nabla\times\tilde{\mbox{\boldmath $V$}}
=\nabla(\nabla\cdot\tilde{\mbox{\boldmath $V$}})-\triangle\tilde{\mbox{\boldmath $V$}}$,
we obtain
$$\begin{array}{c}
\displaystyle
\frac{1}{\mu\epsilon}\nabla\times\nabla\times\tilde{\mbox{\boldmath $V$}}(x)+\tau^2\tilde{\mbox{\boldmath $V$}}(x)\\
\\
\displaystyle
=\mbox{\boldmath $L$}(x,\nabla\mbox{\boldmath $V$}(x^r),\mbox{\boldmath $V$}(x^r))-\frac{2}{\mu\epsilon}d_{\partial D}(x)
\mbox{\boldmath $W$}(x,(\nabla^2\mbox{\boldmath $V$})(x^r),
(\nabla\mbox{\boldmath $V$})(x^r),(\mbox{\boldmath $V$})(x^r)),
\end{array}
\tag {A.20}
$$
where
$$\begin{array}{c}
\displaystyle
\mbox{\boldmath $L$}(x,\nabla\mbox{\boldmath $V$}(x^r),\mbox{\boldmath $V$}(x^r))\\
\\
\displaystyle
=\frac{2}{\mu\epsilon}\left\{(\nabla\cdot\mbox{\boldmath $n$})(x^r)(I-2\pi(x))(\nabla B)(x^r)
+B(x^r)
(I-2\pi(x))(\nabla(\nabla\cdot\mbox{\boldmath $n$}))(x^r)\right\}\\
\\
\displaystyle
-\frac{2}{\mu\epsilon}\mbox{\boldmath $n$}(x)\text{Trace}\,\left(\mbox{\boldmath $V$}'(x^r)\mbox{\boldmath $n$}'(x)\right)\}
\\
\\
\displaystyle
-\frac{2}{\mu\epsilon}\left\{\mbox{\boldmath $n$}'(x^r)(\nabla B)(x^r)+B(x^r)(I-\pi(x))(\triangle\mbox{\boldmath $n$})(x^r)\right\}\\
\\
\displaystyle
-\frac{2}{\mu\epsilon}\left(2\text{Trace}\,(\mbox{\boldmath $A$}'(x^r)\mbox{\boldmath $n$}'(x^r))
+\mbox{\boldmath $A$}(x^r)\cdot(\triangle\mbox{\boldmath $n$})(x^r)\right)
\mbox{\boldmath $n$}(x)\\
\\
\displaystyle
-\frac{2}{\mu\epsilon}(\nabla\cdot\mbox{\boldmath $n$})(x)
\mbox{\boldmath $n$}'(x)\mbox{\boldmath $A$}(x^r)
\\
\\
\displaystyle
+\frac{2}{\mu\epsilon}(\nabla\cdot\mbox{\boldmath $n$})(x)((\nabla B)(x^r)\cdot\mbox{\boldmath $n$}(x))
\mbox{\boldmath $n$}(x)\\
\\
\displaystyle
-\frac{2}{\mu\epsilon}(\mbox{\boldmath $n$}'(x)-
\mbox{\boldmath $n$}'(x^r))(\nabla B)(x^r)-\frac{1}{\mu\epsilon}B(x^r)((\triangle\mbox{\boldmath $n$})(x)
-(\triangle\mbox{\boldmath $n$})(x^r))
\end{array}
$$
and
$$\begin{array}{c}
\displaystyle
\mbox{\boldmath $W$}(x,\nabla^2\mbox{\boldmath $V$}(x^r),
\nabla\mbox{\boldmath $V$}(x^r),\mbox{\boldmath $V$}(x^r))
=\mbox{\boldmath $Z$}(x,(\nabla^2\mbox{\boldmath $V$})(x^r),(\nabla\mbox{\boldmath $V$})(x^r),(\mbox{\boldmath $V$})(x^r))
\\
\\
\displaystyle
+\mbox{\boldmath $T$}(x,(\nabla^2\mbox{\boldmath $A$})(x^r),(\nabla\mbox{\boldmath $A$})(x^r))
-\mbox{\boldmath $N$}(x,(\nabla^2B)(x^r),(\nabla B)(x^r)).
\end{array}
$$

\subsubsection{Computation of $(1/\mu\epsilon)\nabla\times\nabla\times\mbox{\boldmath $C$}+\tau^2\mbox{\boldmath $C$}$
for $\mbox{\boldmath $C$}$ given by (4.7).}

We assume that $\partial D$ is $C^4$.
Set $d=d_{\partial D}(x)$.

Since $\nabla\cdot\mbox{\boldmath $C$}=2d\nabla\cdot (\mbox{\boldmath $n$}'\mbox{\boldmath $A$}(x^r))$,
we have
$$\begin{array}{c}
\displaystyle
\nabla(\nabla\cdot\mbox{\boldmath $C$})
=2\nabla\cdot(\mbox{\boldmath $n$}'\mbox{\boldmath $A$}(x^r))\mbox{\boldmath $n$}
+2d\nabla\{\nabla\cdot(\mbox{\boldmath $n$}'\mbox{\boldmath $A$}(x^r))\}.
\end{array}
\tag {A.21}
$$
On the other hand, we have
$$\begin{array}{c}
\displaystyle
(\triangle\mbox{\boldmath $C$})^i
=2\sum_{j}\triangle\left(d\frac{\partial n^i}{\partial x_j}\right)\mbox{\boldmath $A$}^j(x^r)
+4\sum_{j}\nabla\left(d\frac{\partial n^i}{\partial x_j}\right)\cdot\nabla(\mbox{\boldmath $A$}^j(x^r))
+2\sum_{j}d\frac{\partial n^i}{\partial x_j}\triangle(\mbox{\boldmath $A$}^j(x^r)).
\end{array}
$$
The third term on this right-hand side is the $i$-th component of
$\displaystyle2d\mbox{\boldmath $n$}'(x)\triangle(\mbox{\boldmath $A$}(x^r))$.
However, by (A.14), this is equal to
$$\begin{array}{c}
\displaystyle
2d\mbox{\boldmath $n$}'(x)(\triangle\mbox{\boldmath $A$})(x^r)-4d\mbox{\boldmath $n$}'(x)(\nabla\cdot\mbox{\boldmath $n$})(x)\mbox{\boldmath $A$}(x^r)
-4d^2\mbox{\boldmath $n$}'(x)
\mbox{\boldmath $T$}(x,(\nabla^2\mbox{\boldmath $A$})(x^r),(\nabla\mbox{\boldmath $A$})(x^r)).
\end{array}
$$
Thus, using formula $\nabla\times\nabla\times\mbox{\boldmath $C$}
=\nabla(\nabla\cdot\mbox{\boldmath $C$})-\triangle\mbox{\boldmath $C$}$, we obtain
$$\begin{array}{c}
\displaystyle
\frac{1}{\mu\epsilon}\nabla\times\nabla\times\mbox{\boldmath $C$}+\tau^2\mbox{\boldmath $C$}
=-2
d\mbox{\boldmath $n$}'\left(\frac{1}{\mu\epsilon}(\triangle\mbox{\boldmath $A$})(x^r)-
\tau^2\mbox{\boldmath $A$}(x^r)\right)\\
\\
\displaystyle
+(\text{first and zero-th order terms})
+2d\times(\text{second, first and zero-th order terms}).
\end{array}
$$
The point is the first term on this right-hand side.
From (A.18) and (A.19) we see that
this right-hand side consists of at most
first order terms.

Summing up, we obtain
$$\begin{array}{c}
\displaystyle
\frac{1}{\mu\epsilon}\nabla\times\nabla\times\mbox{\boldmath $C$}+\tau^2\mbox{\boldmath $C$}
=\sum_{ijk}Q_{ijk}(x)\frac{\partial \mbox{\boldmath $A$}^j}{\partial x_k}(x^r)
+\sum_{ij}Q_{ij}(x)\mbox{\boldmath $A$}^j(x^r)\\
\\
\displaystyle
+2d_{\partial D}(x)
\left(\sum_{j,k,l}R_{ijkl}(x)\frac{\partial^2\mbox{\boldmath $V$}^j}{\partial x_k\partial x_l}(x^r)
+\sum_{j,k}R_{ijk}(x)\frac{\partial\mbox{\boldmath $V$}^j}{\partial x_k}(x^r)
+\sum_{j}R_{ij}(x)\mbox{\boldmath $V$}^j(x^r)\right).
\end{array}
\tag {A.22}
$$
Note that all the coefficients are independent of $\tau$ and continuous in a tubular neighbourhood of $\partial D$,
in particular, $R_{ijkl}(x)$ which come from the second order terms in (A.21) is $C^1$.

\vskip1cm
\noindent
e-mail address

ikehata@amath.hiroshima-u.ac.jp

\end{document}